\documentclass[draftcls,onecolumn,journal]{IEEEtran}
\makeatletter
\def\endthebibliography{%
	\def\@noitemerr{\@latex@warning{Empty `thebibliography' environment}}%
	\endlist
}
\usepackage{cite}
\usepackage{amsmath,amssymb,amsfonts}
\usepackage{graphicx}
\usepackage{textcomp}
\usepackage{xcolor}

\usepackage{hyperref}
\usepackage{amsbsy}
\usepackage{amsthm}

\usepackage{algorithmicx}
\usepackage[]{algorithm}
\usepackage{algpseudocode}
\usepackage{comment}

\allowdisplaybreaks
\usepackage{mathtools}
\usepackage{bbm}
\usepackage{xcolor}
\usepackage{footnote}
\providecommand{\abs}[1]{\lvert#1\rvert}
\renewcommand{\b}[1]{\ensuremath{\mathbf{#1}}} 
\newcommand{\bs}[1]{\ensuremath{\boldsymbol{#1}}} 
\renewcommand{\c}[1]{\ensuremath{\mathcal{#1}}} 
\newcommand{\Ex}[1]{\ensuremath{\mathbb{E}[#1]}}  
\newcommand{\ind}{1\hspace{-1.6mm}1} 
\newcommand{\norm}[1]{\ensuremath{\left\|#1\right\|}} 
\newcommand{\hb}[1]{\ensuremath{\hat{\mathbf{#1}}}} 
\newcommand{\eqtext}[1]{\ensuremath{\stackrel{\text{#1}}{=}}} 
\newcommand{\leqtext}[1]{\ensuremath{\stackrel{\text{#1}}{\leq}}} 
\providecommand{\ip}[2]{\langle #1, #2 \rangle} 

\theoremstyle{plain}
\newtheorem{theorem}{Theorem}
\newtheorem{lemma}{Lemma}
\newtheorem*{lemma*}{Lemma}

\newtheorem{cor}{Corollary}

\theoremstyle{definition}

\newtheorem*{defn*}{Definition}

\theoremstyle{remark}
\newtheorem*{rem}{Remark}

\newcommand{\nn}{\nonumber}

\newcommand{\lp}{\left(}
\newcommand{\rp}{\right)}

\newcommand{\lb}{\left[}
\newcommand{\rb}{\right]}

\newcommand{\nf}{\nabla f}
\newcommand{\nfi}{\nabla f_i}
\newcommand{\nuu}{\nabla u}
\newcommand{\nui}{\nabla u_i}
\newcommand{\nfh}{\nabla \fh}

\newcommand{\sumjn}{\sum_{j=1}^n}
\newcommand{\sumin}{\sum_{i=1}^n}
\newcommand{\sumtT}{\sum_{t=1}^T}

\def \a {{\b{a}}}
\def \x {{\b{x}}}
\def \y {{\b{y}}}

\def \z {{\b{z}}}

\def \h {{\b{h}}}

\def \v {{\b{v}}}
\def \w {{\b{w}}}

\def \bo {\mathbf{1}}

\def \tf {{\tilde{f}}}

\def \hx {{\hb{x}}}

\def \A {{\b{A}}}

\def \I {{\b{I}}}
\def \J {{\b{J}}}

\def \W {{\b{W}}}

\def \Wu {{\underline{\W}}}
\def \bou {{\underline{\bo}}}

\def \cE {{\c{E}}}
\def \cH {{\c{H}}}

\def \cN {{\c{N}}}
\def \cV {{\c{V}}}
\def \cX {{\c{X}}}

\def \xib {{\bs{\xi}}}

\def \nab {{\bs{\nabla}}}

\def \delb {{\bs{\delta}}}

\def \EE {{\mathbb{E}}}
\def \Rn {{\mathbb{R}}}

\def \xib {{\boldsymbol{\xi}}}
\def \lamW {{\lambda_{\W}}}

\def \xh {{\hat{\x}}}
\def \fh {{\hat{f}}}

\def \bx {{\bar{\x}}}
\def \by {{\bar{\y}}}
\def \bz {{\bar{\z}}}

\def \bnu {{\bar{\nabla u}}}
\def \bxh {{\bar{\xh}}}

\newtheorem{assumption}{}

\def\BibTeX{{\rm B\kern-.05em{\sc i\kern-.025em b}\kern-.08em
    T\kern-.1667em\lower.7ex\hbox{E}\kern-.125emX}}

\begin{document}
	\title{Analysis of Decentralized Stochastic Successive Convex Approximation for composite non-convex problems
		\author{Basil M. Idrees*, Shivangi Dubey Sharma*\footnote{* both the authors have equal contribution}, and  Ketan Rajawat}}
	\maketitle
	
\maketitle

\begin{abstract}
Successive Convex approximation (SCA) methods have shown to improve the empirical convergence of non-convex optimization problems over proximal gradient-based methods. SCA uses a strongly convex surrogate and offers a more flexible framework to solve such optimization problems. Further, in decentralized optimization, which aims to optimize a global function using only local information, the SCA framework has been successfully applied to achieve improved convergence. Still, the stochastic first-order (SFO) complexity of decentralized SCA algorithms has remained understudied. While non-asymptotic convergence analysis has been studied for decentralized deterministic settings, its stochastic counterpart has only been shown to converge asymptotically.
 
We have analyzed a novel accelerated variant of the decentralized stochastic SCA 
that minimizes the sum of non-convex (possibly smooth) and convex (possibly non-smooth) cost functions. The algorithm viz. \textbf{D}ecentralized \textbf{M}omentum-based \textbf{S}tochastic \textbf{SCA} (\textbf{D-MSSCA}), iteratively solves a series of strongly convex subproblems at each node using one sample at each iteration. 
The recursive momentum-based updates at each node contribute to achieving stochastic first order (SFO) complexity of $\mathcal{O}(\epsilon^{-3/2})$  provided that the step sizes are smaller than the given upper bounds. Even with one sample used at each iteration and a non-adaptive step size, the rate is at par with the SFO complexity of decentralized state-of-the-art gradient-based algorithms. The rate also matches the lower bound for the centralized, unconstrained optimization problems. Through a synthetic example, the applicability of D-MSSCA is demonstrated.
\end{abstract}


\section{Introduction}
 We consider the following decentralized stochastic non-convex composite optimization problem:
\begin{align}  \label{Prob}
  U^\star = \underset{\x \in \Rn^d  }{\min}~~~U(\x)&: =  \frac{1}{n}\sumin u_i(\x) +  h(\x) \tag{$\mathcal{P}$}\\
        \text{s.t.}&~~~g(\x) \leq 0  \nn
\end{align}

where $ u_i(\x) := \Ex{f_i(\x,\xib_i)}$ with $f_i: \Rn^d \rightarrow \Rn$ being smooth but possibly non-convex. On the other hand, $h: \Rn^d \rightarrow \Rn$ is convex but possibly non-smooth while $g: \Rn^d \rightarrow \Rn $ is a convex function. The agents communicate over undirected links, denoted by tuples of the form $(i,j)$. The topology of the communication network is represented by an undirected graph $\mathcal{G} = (\cV, \cE)$, where the node set $\cV$ collects the set of agents and the edge set $\cE$ collects the set of communication links.

	We consider the fully decentralized setting, where $f_i$ is private to agent $i$, though the regularizer $h$ and the constraint function $g$ are public knowledge. Nevertheless, the agents are unaware of the global objective function but must communicate with each other in order to solve \eqref{Prob}. The decentralized stochastic optimization problem in \eqref{Prob} arises in a number of areas, including statistical inference, machine learning, and sensor networks \cite{bottou2018optimization}. In all of these applications, the objective $f_i$ depends on the data or observation $\xib_i$ which is private to agent $i$ and cannot be shared with the other agents. Instead, the agents may only share the derived quantities such as gradients or iterates. 

Existing methods to solve \eqref{Prob} include projected and proximal stochastic gradient methods \cite{bianchi2012convergence, swenson2022distributed, wang2021distributed, xin2021stochastic, yan2023compressed, mancino2023proximal } and Successive Convex Approximation (SCA) \cite{zheng2023distributed, di2019distributed} methods. The performance of these algorithms is measured in terms of the number of stochastic first order (SFO) oracle calls required to reach a $\epsilon$-Karush-Kuhn Tucher (KKT) point. While \cite{bianchi2012convergence} proposes a projected DSGD-type algorithm for problems with a compact constraint set, \cite{swenson2022distributed} goes ahead and establishes the asymptotic convergence of DSGD for a family of non-convex, non-smooth functions. Further, in \cite{wang2021distributed}, a decentralized stochastic proximal primal-dual method called SPPDM is proposed, assuming that the epigraph of $h$ is a polyhedral set. Only three works address non-asymptotic iteration complexity analysis for stochastic non-convex composite problems with a general convex non-differentiable regularizer $h$. While DProxSGT\cite{yan2023compressed} achieves sub-optimal rate of $\mathcal{O}(\epsilon^{-2})$ without mean-squared smoothness assumption, ProxGT-SR-O/E \cite{xin2021stochastic} and DEEPSTORM \cite{mancino2023proximal} achieve an optimal convergence rate of $\mathcal{O}(\epsilon^{-3/2})$. However, the problem with ProxGT-SR-O/E \cite{xin2021stochastic} is that it uses large batches.
 
Unlike the above methods, SCA methods offer a more flexible framework to solve non-convex optimization problems since the first work done by  \cite{scutari2013decomposition}. At each iteration, SCA solves a convexified sub-problem formed by approximating the non-convex functions using convex functions called \textit{surrogates}. Different from other competitive algorithms like Expectation-Minimization (EM) and Majorization-Minimization (MM) SCA offers more freedom in the choice of surrogates which can be tailored to a specific problem at hand \cite{scutari2016parallel, scutari2016paralel}. Even though there is a rich body of work on SCA \cite{yang2016parallel, liu2019stochastic,liu2018online,ye2019stochastic,liu2019two,liu2021two,
 mokhtari2017large,koppel2018parallel, mokhtari2020high, idrees2021practical, idrees2024constrained}, their non-asymptotic analysis has largely remained understudied. Under stochastic centralized settings, through non-asymptotic convergence analysis of SCA it was shown that AsySCA \cite{idrees2021practical} archives a rate of $\mathcal{O}(\epsilon^{-2})$. Further, combining accelerated momentum-based updates with SCA has improved the rate to $\mathcal{O}(\epsilon^{-3/2})$ in \cite{idrees2024constrained}.

Under decentralized settings, there are only a handful of SCA algorithms \cite{zheng2023distributed,di2016next, di2019distributed} including both deterministic and stochastic cases.
Authors in \cite{di2016next, di2019distributed}, proposed an SCA-based decentralized algorithm NEXT, and its stochastic variant S-NEXT, \cite{di2019distributed}.
In both of these works, only asymptotic convergence has been proven. Recently, \cite{zheng2023distributed} introduced a decentralized algorithm that employs Nesterov-like momentum, providing the first non-asymptotic analysis of decentralized SCA methods for deterministic case. However, the rate metric used is not general but SCA-specific, and detailed proofs were only provided for the case where the objective function is convex. To the best of our knowledge, there is no comprehensive non-asymptotic convergence analysis for the decentralized stochastic SCA algorithm in the literature.


	In this work, we have analyzed a novel Decentralized Momentum-based Stochastic SCA (D-MSSCA) algorithm to solve \eqref{Prob}. The D-MSSCA hinges on SCA technique and iteratively solves a convexified subproblem at each node. The recursive momentum type local gradient estimate along with global gradient tracking allowed us to achieve an optimal convergence rate of $\mathcal{O}(\epsilon^{-3/2})$. Our analysis extends the methods used in gradient-based approaches \cite{xin2021stochastic,mancino2023proximal} to SCA framework. One of the challenges in the convergence analysis is obtaining the global descent direction using local or partial information. We overcome this by using the local optimality condition of the strongly convex subproblem at each node. This, in turn, allowed us to develop a recursive bound on the global function $U$ characterizing the progress in a single iteration. Finally, simulations on a synthetic problem empirically validate the theoretical findings.

	A comparative performance of various state-of-the-art algorithms that can be used to solve \eqref{Prob} is provided in Table \ref{litts}. It can be observed that the proposed D-MSSCA algorithm achieves the optimal convergence rate.

\subsection{Notations}
	We denote vectors (matrices) using lowercase (uppercase) bold font letters. The $i$-th entry of vector $\x$ is denoted by $[\x]_i$ while the $(i,j)$-th entry of matrix $\A$ is denoted by $A_{ij}$. The $n\times n$ identity matrix is denoted by $\I_n$ while the $n \times 1$ all-one vector is denoted by $\bo_n$. The Kronecker product is denoted by $\otimes$. The $d$-dimensional average of any $nd$-dimensional vector $\mathbf{a} \in \Rn^{nd}$, is represented by $\bar{\mathbf{a}}=\frac{1}{n}\lp \bo_n^\mathsf{T} \otimes \I_d \rp \mathbf{a} \in \Rn^d $. The Euclidean norm of $\y$ is denoted by $\norm{\y}$. The maximum eigenvalue of a matrix $\A$ is denoted by $\lambda_{\max}(\A)$. The set of subgradients of a function $v:\Rn^d \rightarrow \Rn$ calculated at $\x = \a$ is denoted by $\partial \lp v \rp \mid_{\x = \a}$. For the sake of brevity, we define $\nabla_\x \tf(\x,\x_i^t,\xib_i^t) \mid_{\x = \mathbf{a}} := \nab \tf(\mathbf{a},\x_i^t,\xib_i^t)$. 
 The indicator function is defined as $ \ind_{\cX}(\x)$ where $ \ind_{\cX}(\x) = 0$ if $\x \in \cX$, otherwise $\ind_{\cX}(\x) = \infty$. Finally, $\cH^t$ represents the history of the system generated by $\{\xi_{i}^{\tau}\}_{i=\{1,2,..,n\}}^{\tau \leq t-1}$.

The rest of the paper is organized as follows: Section \ref{propme} discusses the proposed algorithm and the assumptions on the problem are discussed. In Section \ref{conv_ana}, the convergence proof of the proposed algorithm is presented. Section \ref{resultss} briefly describes the proposed algorithm's applicability to a synthetic problem.

 \begin{table*}[]
\centering
\caption{\label{litts} Comparison of oracle complexities of Decentralized consensus stochastic non-convex composite optimization algorithms for Expectation (population risk) problems. (To make comparisons fair, we have converted the SFO-complexities of all the algorithms to match our definition of \eqref{opt})}
	\setlength\tabcolsep{4.5pt}	
\resizebox{\textwidth}{!}{%
\begin{tabular}{|c|c|c|c|}
\hline
\textbf{Algorithm}                                                                                                        & \textbf{SFO complexity}                                                                & \textbf{Asymptotic/ Non-Asymptotic} & \textbf{Remarks}          \\ \hline
projected DSGD \cite{bianchi2012convergence}                                                                                      & -                                                 & Asymptotic             & compact constraint
set \\ \hline
\cite{swenson2022distributed}                                                                                      & -                                                 & Asymptotic             & family of non-convex nonsmooth functions\\ \hline
 SPPDM \cite{wang2021distributed}                                                                                                              & $\mathcal{O}(\epsilon^{-2})$                                                                            & Non-Asymptotic                        & epigraph of $h$ is polyhedral             \\ \hline
 ProxGT-SR-O/E \cite{xin2021stochastic}                                                                                                                & $\mathcal{O}(\epsilon^{-3/2})$                                                                         & Non-Asymptotic                         & Multiple communication, large batches  per iteration             \\ \hline
S NEXT    \cite{di2019distributed}                                                                                                          & -                                                                            & Asymptotic                         & SCA based              \\ \hline
DProxSGT \cite{yan2023compressed}                                                                            & 
$\mathcal{O}(\epsilon^{-2})$ & Non-Asymptotic                    & without MSS assumption            \\ \hline
DEEPSTORM \cite{mancino2023proximal}                                                                            & 
$\mathcal{O}(\epsilon^{-3/2})$ & Non-Asymptotic                    & gradient-based              \\ \hline
\textbf{D-MSSCA (This work)}                                                                           & 
$\mathcal{O}(\epsilon^{-3/2})$ & Non-Asymptotic                    & SCA based           \\ \hline
\end{tabular}%
}
\end{table*}

\section{Proposed method} \label{propme}

\subsection{Problem}

Consider a network of $n$ agents or nodes communicating over a fixed undirected graph $\mathcal{G} = (\cV, \mathcal{E})$, where $\cV$ is the set of nodes and $\mathcal{E}$ is a set of edges or links. An edge $(i, j) \in  \mathcal{E}$ represents a communication link between nodes $i$ and $j$. 
Now we re-write the decentralized problem \eqref{Prob} as,
  \begin{align}  \label{Prob1}
    \underset{\x \in \cX }{\min}~~~  \frac{1}{n}\sumin u_i(\x) +  h(\x) \tag{$\mathcal{P}_c$}
\end{align}  
where $\cX = \{\x \in \Rn^d \mid g(\x) < 0 \} $.

\subsection{Proposed Algorithm}
We will now state the proposed algorithm. All nodes are initialized at the same feasible point $\bx^1 \in \cX $, which satisfies $g(\bx^1) \leq 0 $. Each node $i$ constructs a strong convex surrogate $\fh_i$ and solves the following optimization problem based upon the private knowledge of $f_i$ and public knowledge of $h$ and $g$;
        \begin{align}
         \label{eq:alg_up_eq_1}
            \xh_i^t = \underset{\x_i \in  \mathcal{X}}{\arg \min} & ~~\tf_i \lp \x_i, \x_i^t, \xi_i^t \rp + \ip{\pi_i^t}{\x_i - \x_i^t} + h(\x_i) 
        \end{align} 
        where $\fh$ is strongly convex and 
        \begin{align} \label{def:tf}
            \tf_i \lp \x_i, \x_i^t, \xi_i^t \rp = \fh_i \lp \x_i, \x_i^t, \xi_i^t \rp  + (1 - \beta) \ip {\z_i^{t-1} - \nfi(\x_i^{t-1},\xi_i^{t})}{ \x_i-\x_i^t} 
        \end{align}
        with $\pi_i^t = \y_i^t - \z_i^t$ and  $\beta$ is the step size. One choice of $\fh$ can be $\fh_i \lp \x_i, \x_i^t, \xi_i^t \rp = f_i \lp \x_i^t, \xi_i^t \rp + \nfi \lp \x_i^t, \xi_i^t \rp \lp \x_i-\x_i^t \rp + \frac{\mu}{2} \norm{\x_i - \x_i^t}^2$. Each node then performs the following two updates,
        \begin{align}
             \x_i^{t+1} &  = \sumjn W_{i,j} \v_j^t  =\sumjn W_{i,j} \lp \x_j^t + \alpha \lp \xh_j^t - \x_j^t \rp \rp. \label{eq:alg_up_eq_2} \\
              \z_i^{t+1} &= \nfi(\x_i^{t+1},\xi_i^{t+1}) + (1-\beta) \lp \z_i^{t} - \nfi(\x_i^{t},\xi_i^{t+1})\rp  \label{eq:alg_up_eq_3}
        \end{align}
        This update is inspired from \cite{xin2021hybrid}. It should be noted that in \eqref{eq:alg_up_eq_3} we have used \textit{local} momentum-based gradient estimator $\z_i^t$ \cite{cutkosky2019momentum,tran2019hybrid}. Also, it is noteworthy that the update \eqref{eq:alg_up_eq_3} can also be seen as a convex combination of vanilla SGD  and SARAH-type gradient estimator \cite{nguyen2017sarah}.  Finally, using gradient tracking\cite{sharma2024optimized, qu2017harnessing}, each node updates the local estimate of the global gradient:
        \begin{equation} \label{eq:alg_up_eq_4} 
            \y_i^{t+1} = \sumjn W_{i,j} \lp \y_j^t + \z_j^{t+1} - \z_j^{t} \rp. 
        \end{equation}
The D-MSSCA algorithm is summarised in Algorithm \ref{alg:DMSCA}:
\begin{algorithm}
    \caption{\textbf{D}ecentralized -\textbf{M}omentum based \textbf{S}tochastic \textbf{SCA} (\textbf{D-MSSCA}) at each node $i$ }\label{alg:DMSCA}
    \begin{algorithmic}[1]
        \State \textbf{Require} $
        \x_1^1=\x_2^1= \dots = \x_n^1, \alpha,  \beta>0, \quad \tau_i, \{ w_{ij}\}_{j=1}^n,  \quad   \text{Sample } \xi_i^1,  \z_i^0=\nfi(\x_i^0, \xi_i^1)=0, \y_i^1 =\z_i^1 = \nfi(\x_i^1, \xi_i^1), $
        \For{ $t=1$ \textbf{to} $T$} 
        \State  Minimize local surrogate as per \eqref{eq:alg_up_eq_1} 

        \State Obtain local update of the solution as per \eqref{eq:alg_up_eq_2}
        
        \State Sample  $\xi_i^{t+1}$ and update the local gradient estimates as per \eqref{eq:alg_up_eq_3}  
         
        \State Update the global gradient estimates as per \eqref{eq:alg_up_eq_4}
        \EndFor
        \State \textbf{Output} $\tilde{\x}_T$ selected uniformly at random from $\{\xh_i^t\}_{0\leq t \leq T}^{i \in \cV}$
        
    \end{algorithmic}
\end{algorithm}

\subsection{Assumptions}
We will now state the assumptions required for the proposed algorithm. The assumptions are divided among the following 3 heads, viz assumption on \eqref{Prob}, those on the surrogate, and those on the network,
\subsubsection{Assumptions on \eqref{Prob}}
\begin{assumption} \label{assm:below}
    The objective $U$ is bounded below, i.e., $\underset{\x \in \Rn^d}{\inf} U(x)> -\infty$
\end{assumption} 
\begin{assumption} \label{assm:exp_equal_true}
Let $\cH^t$ represent the history of the system generated by $\{\xi_{i}^{\tau}\}_{i=\{1,2,..,n\}}^{\tau \leq t-1}$, then 
    $$\EE\lb \nfi(\x^{t},\xi_i^{t} \mid \cH^t ) \rb = \nui(\x^{t}); $$ 
\end{assumption}
\begin{assumption} \label{assm:bounded_var}
Bounded Variance:
    $\EE\lb \norm{\nfi(\x,\xi_i^{t}) - \nui(\x)}^2 \rb \leq \sigma_i^2~~~ \forall \x \in \Rn^d, \bar{\sigma}^2 = \sumin \sigma_i^2 $
\end{assumption}

\begin{assumption} \label{assm:smoothness}
Mean Squared Smoothness(MSS): Each local function $f_i$ is $L$-smooth ,
    $$\EE \norm{\nf_i(\x^t,\xi^t)-\nf_i(\y^t,\xi^t)} = L \EE \norm{\x^t - \y^t}; $$ 
\end{assumption}

Assumptions \ref{assm:below}- \ref{assm:smoothness} are standard in the context of distributed optimization. A direct consequence of \ref{assm:below} is that for an initial point $\x_i^1 \in \cX $  we have $U (\bx^1) - U^\star \leq B $. Assumption \ref{assm:smoothness} implies that $u$ is also $L$-smooth. is introduced to simplify the analysis.
\subsubsection{Assumptions on the surrogate} 
Two assumptions on the surrogate choice are:
\begin{assumption} \label{assm:tangent}
    Tangent matching: $\nfh_i(\x^t,\x^{t},\xi_i^{t}) = \nfi(\x^{t},\xi_i^{t}); $ 
\end{assumption}
\begin{assumption}\label{assm:muconvex}
Each surrogate $\fh_i$ of local function $f_i$  is $\mu-$ strongly convex.
\end{assumption}
Assumptions \ref{assm:tangent} and  \ref{assm:muconvex} are standard in the context of SCA and restrict the choice of surrogates. A consequence of \ref{assm:tangent} is that $\nabla \tf_i \lp \x_i^t, \x_i^t, \xi_i^t \rp = \z_i^t$.
\subsubsection{Assumption on Network}
\begin{assumption}\label{assm:condn_on_graph}
Graph $\mathcal{G} \lp \cV, \mathcal{E} \rp$ is undirected, connected and communication matrix $\W$ is doubly stochastic. $W_{i,i}>0$ for all $i$ in $\cN$ and for $i \neq j$, $w_{ij} > 0$ wherever $\lp i, j \rp \in \mathcal{E}$, $W_{i,j} = 0$ otherwise.
\end{assumption}
Assumption \ref{assm:condn_on_graph} is standard in the context of decentralized optimization and is required to achieve consensus among agents.

\subsection{Approximate optimality}
The performance of the proposed algorithm is studied in terms of its SFO complexity. We define the following metric, viz \emph{mean-squared stationary gap}\cite{xin2021improved}, that provides the number of calls to the SFO oracle to achieve an $\epsilon$-KKT point in expectation.
\begin{align}
     \frac{1}{n} \sumin &\frac{1}{T} \sum_{t=1}^{T} \EE \lb \norm{\nabla u(\xh_i^t) +\hat{\w}_i^t }^2 \rb \leq \epsilon \label{opt}
\end{align}
where $ \hat{\w}_i^t \in \partial (h + \ind_{\cX}) \mid_{\x =\xh_i^t}$. If \eqref{opt} holds then the output $\tilde{\x}$ of D-MSSCA is chosen uniformly at random from the set  $\{\xh_i^t\}_{0\leq t \leq T}^{i \in \cV}$ then we have $\EE [ \|\nabla u(\tilde{\x}) + \hat{\w}_{\tilde{\x}}\|^2 ] \leq \epsilon$ where $ \hat{\w}_{\tilde{\x}} \in \partial (h + \ind_{\cX}) \mid_{\x = \tilde{\x}}$.

\section{Convergence Analysis}  \label{conv_ana}
This section will provide a detailed convergence analysis of the D-MSSCA algorithm and compare its rate with other state-of-the-art algorithms. For the analysis, we define the concatenated vectors $\x, \y, \z, \xh \in \Rn^{nd}$ by concatenating corresponding local vectors $\{\x_i \in \Rn^d\}_{i \in \cV}, \{\y_i \in \Rn^d\}_{i \in \cV}, \{\z_i \in \Rn^d\}_{i \in \cV}, \{\xh_i \in \Rn^d\}_{i \in \cV}$ of all the nodes. Using these concatenated vectors, we can write the update equation \eqref{eq:alg_up_eq_2} and \eqref{eq:alg_up_eq_4} in a more compact form as below,
\begin{align}
     \x^{t+1} &= \Wu \lp \x^t + \alpha \lp \xh^t - \x^t \rp \rp \label{eq:alg_up_cmpct_eq_2}, \\
    \y^{t+1} &= \Wu \lp \y^t + \z^{t+1} - \z^t \rp, \label{eq:alg_up_cmpct_eq_4}
\end{align}
where, $\Wu = \W \otimes \I_d \in \Rn^{nd \times nd}.$ Furthermore, we define the concatenated local gradient vector $\nuu(\x^t) := [\nuu_1(\x_1^t)^\mathsf{T}, \nuu_2(\x_1^t)^\mathsf{T},  \cdots, \nuu_n(\x_n^t)^\mathsf{T}]^\mathsf{T} \in \Rn^{nd}$, where $\nuu_i(\x_i^t) \in \Rn^d $ for all $ i \in \cV$. 
For the sake of brevity, we define for all $t$
\begin{align}
    \theta^t &= \norm{\x^t - \frac{1}{n}\lp \bo_n \otimes \I_d \rp \bx^t}  &\text{(consensus error)} \label{eq:def_consensus_error}, \\
    \delb^t &= \xh^t - \x^t  &\text{(iterate progress)} \label{eq:def_itr_prgrs},\\
    \phi^t &= \EE \lb \norm{\bz^t - \bnu(\x^t)}^2 \rb &\text{(global gradient variance)} \label{eq:defglobal_grad_var},\\
    \upsilon^t &= \EE \lb \norm{\z^t - \nuu(\x^t)}^2 \rb &\text{(network gradient variance)} \label{eq:def_netwrk_grad_var},\\
    \varepsilon^t&= \EE \lb \norm{\y^t - \frac{1}{n}\lp \bo_n \otimes \I_d \rp \by^t}^2 \rb &\text{(gradient tracking error)} \label{eq:def_track_var}.
\end{align}
We begin our analysis by stating some standard results used in decentralized optimization in Lemma~\ref{lem:basic_decentralized_result} \cite{xin2021hybrid,xin2021stochastic,xin2022fast}. Next, Lemmm~\ref{lem:consnss_bound}, Lemma~\ref{lem: error_in_grad_esti} and Lemma~\ref{lem:bndnss_y} establish the contraction relationships for $\theta^t, \phi^t, \upsilon^t$ and $\varepsilon^t$. These contraction relationships help in bounding the cumulative error accumulation $\sumtT \EE\lb \theta^t\rb, \sumtT \phi^t, \sumtT \varepsilon^t$ and $\sumtT \upsilon^t$ in Lemma~\ref{lem:sum_cnsnss_cmutatv}, Lemma~\ref{lem:HSGD}, and Lemma~\ref{lem:sumy}. Lemma~\ref{lem:HSGD}, then extends these findings to obtain bounds on the accumulated average iterate progress $\sumtT \EE \lb \delta^t \rb$ under certain step-size conditions. Finally, Theorem~\ref{the1} uses all these results to quantify the SFO complexity of D-MSSCA.

The results presented in Lemma~\ref{lem:consnss_bound}-\ref{lem:sumy} can be obtained by applying the approach of \cite{xin2021hybrid} to the SCA framework. However, our bounds are different as we define the quantities in terms of $\delb$ rather than $\upsilon$, and the variations in the intermediate steps are detailed in the proof. The results in Lemma~\ref{lemU} and Lemma~\ref{lem:cmu} are different from \cite{xin2021hybrid} due to our focus on average progress $\Delta^T = \frac{1}{T} \sumtT \EE \lb \norm{\delb^t}^2 \rb$.  The results of Lemma~\ref{lem:cmu} and Theorem~\ref{the1} are entirely novel to this work. It is worth noting that the bound in Lemma~\ref{lem:cmu} is similar to the deterministic bound in \cite{zheng2023distributed}. However, the steps needed to obtain those bounds have been omitted in \cite{zheng2023distributed}. 
 
\begin{lemma} \label{lem:basic_decentralized_result} Under Assumptions  \ref{assm:smoothness} and \ref{assm:condn_on_graph}, we have the following results for all $t\geq1$, where $\x^t, \y^t, \z^t$ are variables of D-MSSCA at iterate $t$.
    \begin{subequations}
    \begin{align}
        \norm{\Wu\x - \frac{1}{n}\lp \bo_n\bo_n^\mathsf{T} \otimes \I_d \rp \x} &\leq \lamW  \norm{\x - \frac{1}{n}\lp \bo_n\bo_n^\mathsf{T} \otimes \I_d \rp \x} \quad \forall \x \in \Rn^{nd}, \label{eq:sub:lem:basic_decentralized_result_1}\\
        \norm{ \sum_{i=1}^n \nuu_i(\bx^t) -  \frac{1}{n}\lp \bo_n^\mathsf{T} \otimes \I_d \rp \nuu(\x^t)}^2 &\leq \frac{L^2 }{n} \norm{\x^t - \frac{1}{n}\lp \bo_n\bo_n^\mathsf{T} \otimes \I_d \rp \x^t}^2, \label{eq:sub:lem:basic_decentralized_result_2} \\
        \by^t &= \bz^t, \label{eq:sub:lem:basic_decentralized_result_3}\\
        \norm{\bx - \by}^2 &\leq \frac{1}{n} \norm{\x - \y}^2 \quad \text{for any} \quad \x,\y \in \Rn^{nd}, \label{eq:sub:lem:basic_decentralized_result_4}
    \end{align}
    \end{subequations}
\end{lemma}
where $\lamW := \lambda_{\max}(\W - \frac{1}{n} \bo \bo^{\mathsf{T}})$.
The proofs of the above results are straightforward and can be found in \cite{scutari2016paralel,qu2017harnessing}. In the proof of \eqref{eq:sub:lem:basic_decentralized_result_1}, the contraction property of doubly stochastic symmetric matrices is used. In \eqref{eq:sub:lem:basic_decentralized_result_2}, the smoothness of the functions $u_i$ ($i \in \cV$) is applied. In \eqref{eq:sub:lem:basic_decentralized_result_3}, the special initialization condition of $v^t$ and the doubly stochastic property of $\W$ are used. Lastly, in  \eqref{eq:sub:lem:basic_decentralized_result_4}, the norm property along with the Cauchy-Schwarz inequality is applied. The next Lemma bounds the consensus errors in the $\x^t-$updates \eqref{eq:alg_up_cmpct_eq_2} of the D-MSSCA algorithm.

\begin{lemma} \label{lem:consnss_bound}
    Under \eqref{assm:condn_on_graph}, for the $\x^t-$ updates of D-MSSCA algorithm, the following inequality holds for all $t \geq 2$ and  $ \eta_1 >  0$
    \begin{align}
       (\theta^t)^2 \leq \lp 1 + \eta_1\rp  \lamW^2 \lp \theta^{t-1}\rp^2 + \lp 1 + \frac{1}{\eta_1} \rp \alpha^2 \lamW^2 \norm{\delb^{t-1}}^2. \nn
    \end{align}
\end{lemma}
The proof of Lemma~\ref{lem:consnss_bound} is provided in Appendix~\ref{pr:lem:consnss_bound} and follows by applying update step~\eqref{eq:alg_up_cmpct_eq_2}, then separating the terms using Young's inequality. Finally, by using the properties of the communication matrix $\W$ \eqref{eq:sub:lem:basic_decentralized_result_1}; we obtain the desired results. Similar contraction bounds on $\theta^t$ have been achieved in various gradient-tracking-based decentralized optimization algorithms \cite{qu2017harnessing,xin2021hybrid}. However, our bound is slightly different because we define it in terms of $\theta^{t-1}$ and $\norm{\delb^{t-1}}^2$ instead of $\theta^{t-1}$ and $\upsilon^{t-1}$, as seen in the literature. Next, we will bound the gradient variances.

\begin{lemma} \label{lem: error_in_grad_esti}
Under \ref{assm:exp_equal_true}-\ref{assm:smoothness}, and \ref{assm:condn_on_graph}, 
The following inequalities hold for the iterates produced by D-MSSCA algorithm, where $t\geq2, 0 < \alpha < 1, \eta_1,\eta_2,\eta_3 > 0$ 
\begin{align} \label{eq:lem:3re3}
        \phi^t  &\leq (1-\beta)^2 \phi^{t-1}  +\frac{3L^2(1- \beta)^2}{n^2} \lp 1 + \frac{1}{\eta_2} \rp \EE \lb \lp \theta^t \rp^2 + n \alpha^2 \norm{\lp \frac{1}{n} \bo^\mathsf{T} \otimes \I_d\rp \delb^{t-1}}^2 + \lp \theta^{t-1}\rp^2 \rb \nn \\
   &+ \frac{(1 +\eta_2)\beta^2 \bar{\sigma}^2}{n^2}, 
\end{align}
and
\begin{align}  \label{eq:lem:3re4}
    \upsilon^t &\leq  (1-\beta)^2 \upsilon^{t-1} +3L^2(1- \beta)^2 \lp 1 + \frac{1}{\eta_3} \rp \EE \lb \lp \theta^t \rp^2 + n \alpha^2 \norm{\lp \frac{1}{n} \bo^\mathsf{T} \otimes \I_d\rp \delb^{t-1}}^2 + \lp \theta^{t-1}\rp^2 \rb \nn \\
 &+ (1 + \eta_3) \beta^2 \bar{\sigma}^2. 
\end{align}
\end{lemma}
The proof of Lemma~\ref{lem: error_in_grad_esti} is identical to that of \cite[Lemma 3]{xin2021hybrid}. However, our final bounds are in terms of $\bar{\delb}^{t-1}$ rather than $\bz^{t-1}$ in \cite{xin2021hybrid}. This difference is due to the $\x-$update of D-MSSCA \eqref{eq:alg_up_eq_2}, which differs from GT-HSGD algorithm proposed in \cite{xin2021hybrid}. The proof uses the unbiased nature of the local gradient estimate $\z_i^t$. Furthermore, as the gradient estimate at each node is independent of those at other nodes given the history sequence $\cH^t$, we can omit the cross terms of inner products appearing in the intermediate steps to obtain simplified expressions. Finally, by applying \ref{assm:smoothness} alongside $\x_i^t$ updates, we get the desired result.

\begin{lemma} \label{lem:bndnss_y}
    Under \ref{assm:exp_equal_true}-\ref{assm:smoothness} and \ref{assm:condn_on_graph}, the following inequality holds for $\beta \in (0,1)$, $\forall t \geq 2$,
    \begin{align}
        \varepsilon^t
        &\leq \frac{1 +\lamW^2}{2} \varepsilon^{t-1} + \frac{4 \beta^2  \lamW^2  }{1- \lamW^2} \upsilon^{t-1} + 3 \lamW^2 \beta^2 \bar{\sigma}^2 + \frac{36 \lamW^2 L^2}{1 - \lamW^2} \EE \lb (\theta^{t-1})^2 \rb + \frac{36 \alpha^2 \lamW^2 L^2 }{1 - \lamW^2} \EE \lb \norm{ \delb^{t-1}}^2 \rb \nn. 
    \end{align}
\end{lemma}
The above lemma bounds the error in local estimation of the global gradient. The proof uses the $\y^t-$update \eqref{eq:alg_up_cmpct_eq_4}, applies conditional expectation, and simplifies intermediate steps using Assumptions~\ref{assm:exp_equal_true},\ref{assm:bounded_var}, and \ref{assm:smoothness}. Finally, by using the consensus error bound in Lemma~\ref{lem:consnss_bound} and performing a few mathematical simplifications, the desired result is achieved. The elaborated proof can be found in Appendix~\ref{pr:lem:bndnss_y}.
The next Lemma provides the basic results of non-negative sequences, which are necessary to bound the error accumulation in subsequent lemmas.

\begin{lemma} \label{lem:seq_bnd} 
The recursions of well-defined sequences can be bound as below: 
    \begin{enumerate} 
        \item Let $\{V^t\}_{t\geq0}$, $\{Q^t\}_{t\geq0}$ be non-negative sequences and $C>0$ be some constant such that $V^t \leq q V^{t-1} + Q^{t-1} + C$ for some $q \in \lp 0,1 \rp$ and for all $t \geq 1$. Then the following inequality holds $\forall T \geq 1$
        \begin{equation}
            \sum_{t=0}^{T} V^t \leq \frac{V^0}{1-q} + \frac{\sum_{t=0}^{T-1} Q^t}{1-q} + \frac{CT}{1-q}. \label{eq:lem:seq_bnd}
        \end{equation}
        \item Let $\{V^t\}_{t\geq1}$, $\{Q^t\}_{t\geq1}$ be non-negative sequences and $C>0$ be some constant such that $V^t \leq q V^{t-1} + Q^{t-1} + C$ for some $q \in \lp 0,1 \rp$ and for all $t \geq 2$. Then the following inequality holds $\forall T \geq 2$
        \begin{equation}
            \sum_{t=1}^{T} V^t \leq \frac{V^1}{1-q} + \frac{\sum_{t=2}^T Q^{t-1}}{1-q} + \frac{CT}{1-q}. \label{eq:lem:seq_bnd_frm1}
        \end{equation}
    \end{enumerate}
\end{lemma}
The above mentioned recursion results align with \cite[Lemma 6]{xin2021hybrid}, and the proof follows a similar approach. However, there is a slight difference between our results and theirs, due to the variation in the relationship of sequences and the range of $T$ considered here. For the complete proof, refer to \cite{xin2021hybrid}. Using the stated lemmas, upper bounds for the cumulative errors up to iteration $T$ can now be established, as detailed in the following lemmas.

\begin{lemma} \label{lem:sum_cnsnss_cmutatv}
    For the proposed D-MSSCA algorithm, following inequality holds: $\forall T > 1$, $\alpha \in \lp 0, 1 \rp$,  
    \begin{align}
        \sum_{t=1}^{T} \EE \lb \lp \theta^t \rp^2 \rb &\leq \frac{4 \alpha^2 \lamW^2 }{\lp 1 -  \lamW^2 \rp^2} \sum_{t=1}^{T-1} \EE \lb \norm{\delb^t}^2 \rb. \nn
    \end{align}
\end{lemma}
The above result can be obtained by summing both sides of the bound obtained in Lemma~\ref{lem:consnss_bound}, for $1 \leq t \leq T$ and applying \eqref{eq:lem:seq_bnd_frm1}. The detailed proof is provided in Appendix~\ref{pr:lem:sum_cnsnss_cmutatv}.

\begin{lemma} \label{lem:HSGD}
For the proposed D-MSSCA algorithm, following inequality holds: $\forall T > 1$, $\beta, \alpha \in (0,1)$,
    \begin{align}
        \sum_{t=1}^T \phi^t 
        & \leq \frac{\bar{\sigma}^2}{n^2 b_0 \beta} + \frac{2 \beta  \bar{\sigma}^2 T}{n^2} +  \frac{6L^2\alpha^2}{n^2\beta}  \sum_{t=1}^{T-1}  \EE \lb \norm{\delb^t}^2 \rb 
        + \frac{12 L^2}{\beta n^2} \sum_{t=1}^{T} \EE \lb \lp \theta^t \rp^2 \rb, \nn
    \end{align}
    \begin{align}
        \sum_{t=1}^T \upsilon^t 
        & \leq \frac{\bar{\sigma}^2}{b_0 \beta}  + 2 \beta \bar{\sigma}^2T    + \frac{6L^2 \alpha^2 }{\beta} \sum_{t=1}^{T-1} \EE \lb \norm{\delb^t}^2  \rb  + \frac{12L^2 }{\beta} \sum_{t=1}^{T} \EE \lb   \lp \theta^t \rp^2 \rb. \nn
    \end{align}
\end{lemma}
 To prove Lemma~\ref{lem:HSGD}, we have applied the results of Lemma~\ref{lem: error_in_grad_esti}, Lemma~\ref{lem:seq_bnd} and \eqref{assm:bounded_var}. The proof of Lemma~\ref{lem:HSGD} is provided in Appendix~\ref{pr:lem:HSGD}.

\begin{lemma} \label{lem:sumy}
    The following inequality holds for all $t>0$
    \begin{align}
         \sum_{t=1}^{T} \varepsilon^t &\leq  \frac{24 \alpha^2 \lamW^2 L^2}{\lp 1- \lamW^2 \rp^2}\lp 3 + 2 \beta \rp \sum_{t=1}^{T-1} \EE \lb \norm{\delb^t}^2  \rb + \frac{24 \lamW^2 L^2}{\lp 1- \lamW^2 \rp^2} \lp 3 + 4 \beta  \rp \sum_{t=1}^{T} \EE \lb   \lp \theta^t \rp^2 \rb  \nn \\
        &+\frac{2 \lamW^2 \beta^2 \bar{\sigma}^2}{1- \lamW^2} \lp \frac{4}{\lp 1- \lamW^2 \rp b_0 \beta} + \frac{8 \beta T }{\lp 1- \lamW^2 \rp}  + 3 T \rp +\frac{2}{1- \lamW^2}  \varepsilon^{1} . \nn
    \end{align}
\end{lemma}
To prove Lemma~\ref{lem:sumy}, we have applied the results of Lemma~\ref{lem:bndnss_y},\ref{lem:seq_bnd} and \ref{lem:HSGD}. The proof of Lemma~\ref{lem:sumy} is provided in Appendix~\ref{pr:lem:sumy}.

The next Lemma is a key result in the analysis of D-MSSCA algorithm, offering a descent inequality for average (over the network) D-MSSCA updates with respect to the global function $U$.

\begin{lemma} \label{lemU}
    Under Assumptions~\ref{assm:smoothness} and \ref{assm:condn_on_graph}, the following inequality holds for all $T > 1$, $\gamma_1 > 0$, $\mu>0$ and $0< \alpha <1$, where $\x^t, \y^t, \z^t$ are iterate variables of D-MSSCA at iterate $t$,
    \begin{align}
       U( \bx^{T+1}) - U(\bx^1)  
       &\leq   \frac{ 3L^2 \alpha  \gamma_1}{2 n}   \sumtT \EE \lb (\theta^t)^2 \rb  + \frac{3 \alpha  \gamma_1}{2} \sumtT \phi^t  + \frac{\alpha}{n} \lp - \mu +\frac{1}{2 \gamma_1 } + \frac{\alpha L}{2} \rp \sumtT \EE \lb \norm{ \delb^t }^2 \rb \nn \\
       &+  \frac{3 \alpha  \gamma_1}{2 n} \sumtT \varepsilon^t. \nn
    \end{align}
\end{lemma} 
The proof of Lemma~\ref{lemU} begins by defining the optimality condition of \eqref{eq:alg_up_eq_1}, and using it with the properties of surrogate to get the descent direction. Further by using the convexity of $\h$ and the properties of the communication matrix $\W$, along with some mathematical simplifications, we achieve the desired result. The complete proof of Lemma~\ref{lemU} is detailed in Appendix~\ref{pr:lemU}.

Now, we will use Lemma~\ref{lem:sum_cnsnss_cmutatv}-\ref{lemU} to upper bound the average progress $\Delta^T = \frac{1}{T} \sumtT \EE \lb \norm{\delb^t}^2 \rb$.

\begin{lemma} \label{lem:cmu} Under considered assumption \ref{assm:below}-\ref{assm:condn_on_graph}, if $ 0<\beta=\alpha^2<1, \mu \geq \frac{6 \sqrt{3} L}{n} \lp 1 + \frac{8 \lamW^2}{(1 - \lamW^2)}\rp$ and  $0 <\alpha \leq \min\Bigg\{\frac{1}{114}, \frac{\lp 1 -  \lamW^2 \rp^2}{432 \lamW^2}, \frac{\lp 1 -  \lamW^2 \rp^{2/3}}{8 \lamW^{2/3}}, \frac{\mu}{6 L} , \frac{\mu^2 \lp 1 -  \lamW^2 \rp^2}{48 L^2 \lamW^2} \Bigg\}  $ then the average progress of D-MSSCA algorithm is upper bounded for all $T\geq 2$ as below
   \begin{align}
        \Delta^T &\leq \frac{4n}{\alpha T \mu}U(\bx^1)- \frac{4n}{\alpha T \mu}U^\star  + \frac{24}{T \mu^2 \lp 1- \lamW^2 \rp} \norm{\nabla u(\x^1)}^2  + \frac{12 \bar{\sigma}^2}{2 T \mu^2} \Bigg(  \frac{1}{b_0 \alpha^2 n} + \frac{2 \alpha^2 T}{n} + \frac{4}{b_0^2 \lp 1- \lamW^2 \rp}\nn \\
        &+ \frac{2 \lamW^2 \alpha^4}{1- \lamW^2} \lp \frac{4}{\lp 1- \lamW^2 \rp b_0 \alpha^2} + \frac{8 \alpha^2 T }{\lp 1- \lamW^2 \rp}  + 3 T \rp \Bigg). \nn
    \end{align}
\end{lemma}
\begin{IEEEproof}
We begin by substituting the bound of $\sum_{t=1}^T \phi^t$ from Lemma~\ref{lem:HSGD} into Lemma~\ref{lemU} and obtain,
    \begin{align}
       U(\bx^{T+1}) - U(\bx^1)  
       &\leq   \frac{ 3L^2 \alpha  \gamma_1}{2 n}   \sumtT \EE \lb (\theta^t)^2 \rb  + \frac{\alpha}{n} \lp - \mu +\frac{1}{2 \gamma_1 } + \frac{\alpha L}{2} \rp \sumtT \EE \lb \norm{ \delb^t }^2 \rb \nn \\
       &+\frac{3 \alpha  \gamma_1}{2 n} \sumtT \varepsilon^t + \frac{3 \alpha  \gamma_1}{2} \lp \frac{\bar{\sigma}^2}{n^2 b_0 \beta} + \frac{2 \beta  \bar{\sigma}^2 T}{n^2} +  \frac{6L^2\alpha^2}{n^2\beta}  \sum_{t=1}^{T-1}  \EE \lb \norm{\delb^t}^2 \rb 
        + \frac{12 L^2}{\beta n^2} \sum_{t=1}^{T} \EE \lb \lp \theta^t \rp^2 \rb \rp. \nn
    \end{align}    
    Combining the common terms and substituting the bound of $\sumtT \varepsilon^t$ from Lemma~\ref{lem:sumy}, we get
    \begin{align}
       U(\bx^{T+1}) - &U(\bx^1) 
       \leq \frac{3 L^2 \alpha \gamma_1}{2n} \lp 1 +  \frac{12}{\beta n}\rp \sumtT \EE \lb (\theta^t)^2 \rb + \frac{\alpha}{n} \lp - \mu +\frac{1}{2 \gamma_1 } + \frac{\alpha L}{2} + \frac{9L^2\alpha^2 \gamma_1}{n \beta} \rp \sumtT \EE \lb \norm{ \delb^t }^2 \rb \nn \\
       &\qquad + \frac{3 \alpha  \gamma_1 \bar{\sigma}^2}{2 n^2} \lp \frac{1}{b_0 \beta} + 2 \beta T \rp +\frac{3 \alpha  \gamma_1}{2 n} \Bigg[ \frac{24 \alpha^2 \lamW^2 L^2}{\lp 1- \lamW^2 \rp^2}\lp 3 + 2 \beta \rp \sum_{t=1}^{T-1} \EE \lb \norm{\delb^t}^2  \rb +\frac{2}{1- \lamW^2}  \varepsilon^{1}  \nn \\
        &\qquad + \frac{24 \lp 3 + 4 \beta  \rp \lamW^2 L^2}{\lp 1- \lamW^2 \rp^2}  \sum_{t=1}^{T} \EE \lb   \lp \theta^t \rp^2 \rb  + \frac{2 \lamW^2 \beta^2 \bar{\sigma}^2}{1- \lamW^2} \lp \frac{4}{\lp 1- \lamW^2 \rp b_0 \beta} + \frac{8 \beta T }{\lp 1- \lamW^2 \rp}  + 3 T \rp \Bigg].  \nn
    \end{align}   
     By combining the common terms and substituting the bound of $\sumtT \EE \lb \lp \theta^t \rp^2\rb$ from Lemma~\ref{lem:sum_cnsnss_cmutatv}, we obtain
    \begin{align}
        U(\bx^{T+1}) - U(\bx^1) &\leq \frac{\alpha}{n} \Bigg[ \frac{6 \alpha^2 \lamW^2  \gamma_1 L^2}{ \lp 1 -  \lamW^2 \rp^2}  \lp  19 + 12 \beta + \frac{12}{\beta n} + \frac{24 \lamW^2 }{\lp 1- \lamW^2 \rp^2} \lp 3 + 4 \beta  \rp \rp  - \mu +\frac{1}{2 \gamma_1 } + \frac{\alpha L}{2} \nn \\
        &+ \frac{9L^2\alpha^2 \gamma_1}{n \beta}  \Bigg] \sumtT \EE \lb \norm{ \delb^t }^2 \rb + \frac{3 \alpha  \gamma_1}{n \lp 1- \lamW^2 \rp} \varepsilon^{1} + \frac{3 \alpha  \gamma_1 \bar{\sigma}^2}{2 n} \Bigg(  \frac{1}{b_0 \beta n} + \frac{2 \beta T}{n} \nn \\
       &+ \frac{2 \lamW^2 \beta^2}{1- \lamW^2} \lp \frac{4}{\lp 1- \lamW^2 \rp b_0 \beta} + \frac{8 \beta T }{\lp 1- \lamW^2 \rp}  + 3 T \rp \Bigg). \nn
    \end{align} 
    Defining $ C_{\mu} =-\Bigg[ \frac{6 \alpha^2 \lamW^2  \gamma_1 L^2}{ \lp 1 -  \lamW^2 \rp^2}  \lp  19 + 12 \beta + \frac{12}{\beta n} + \frac{24 \lamW^2 }{\lp 1- \lamW^2 \rp^2} \lp 3 + 4 \beta  \rp \rp  -\mu +\frac{1}{2 \gamma_1 }  + \frac{\alpha L}{2} + \frac{9L^2\alpha^2 \gamma_1}{n \beta}  \Bigg]$, we have:
    \begin{align}
        \frac{C_{\mu}}{T} \sumtT \EE \lb \norm{ \delb^t }^2 \rb  &\leq \frac{n}{\alpha T}U(\bx^1)- \frac{n}{\alpha T}U(\bx^{T+1}) + \frac{3 \gamma_1}{T \lp 1- \lamW^2 \rp} \varepsilon^{1} + \frac{3 \gamma_1 \bar{\sigma}^2}{2 T} \Bigg(  \frac{1}{b_0 \beta n} + \frac{2 \beta T}{n} \nn \\
       &~~~+ \frac{2 \lamW^2 \beta^2}{1- \lamW^2} \lp \frac{4}{\lp 1- \lamW^2 \rp b_0 \beta} + \frac{8 \beta T }{\lp 1- \lamW^2 \rp}  + 3 T \rp \Bigg). \nn
    \end{align} 
    Also, from the initialization of $\z_i^1$ and the update \eqref{eq:alg_up_eq_4} of the D-MSSCA algorithm, we have:
    \begin{align}
        \varepsilon^1 &= \EE \norm{(\I - \frac{1}{2}\bo \bo^{\mathsf{T}}\otimes \I_d)\y^1}^2 \leq \EE \norm{\y^1}^2 = \EE \norm{\nf (\x^1, \xib^1) - \nuu (\x^1) + \nuu (\x^1)}^2 \nn \\
        &\leq  2\sumin \EE \norm{\frac{1}{b_0} \sum_{r=1}^{b_0} \lp \nf_i (\x_i^1, \xib_i^{1,r}) - \nuu_i (\x_i^1)\rp}^2 + 2\EE \norm{\nuu (\x^1)}^2, \nn \\
        &\leqtext{(i)} \frac{2}{b_0^2} \sumin \sum_{r=1}^{b_0} \EE \norm{ \nf_i (\x_i^1, \xib_i^{1,r}) - \nuu_i (\x_i^1)}^2 + 2\EE \norm{\nuu (\x^1)}^2, \nn \\
        &\leqtext{(ii)} \frac{2 \bar{\sigma}^2}{b_0^2} + 2 \norm{\nabla u(\x^1)}^2, \label{eq:lem:ref6}
    \end{align}
    where in (i), we used the fact that $\xib_i^{1,l}, \xib_j^{1,m}$ are independent for $l \neq m$ and in (ii), we applied \eqref{assm:bounded_var}. Substituting the bound $\varepsilon^1 \leq \frac{2 \bar{\sigma}^2}{b_0^2} + 2 \norm{\nabla u(\x^1)}^2,$ and $U^\star < U(\bx^{T+1}) $, we obtain:, 
    \begin{align}
        \frac{C_{\mu}}{T} \sumtT \EE \lb \norm{ \delb^t }^2 \rb  &\leq \frac{n}{\alpha T}U(\bx^1)- \frac{n}{\alpha T}U^\star  + \frac{6 \gamma_1}{T \lp 1- \lamW^2 \rp} \norm{\nabla u(\x^1)}^2  + \frac{3 \gamma_1 \bar{\sigma}^2}{2 T} \Bigg(  \frac{1}{b_0 \beta n} + \frac{2 \beta T}{n} + \frac{4}{b_0^2 \lp 1- \lamW^2 \rp} \nn \\
       &~~~+ \frac{2 \lamW^2 \beta^2}{1- \lamW^2} \lp \frac{4}{\lp 1- \lamW^2 \rp b_0 \beta} + \frac{8 \beta T }{\lp 1- \lamW^2 \rp}  + 3 T \rp \Bigg). \nn
    \end{align} 
    Further if we consider $\beta = \alpha^2$,$\gamma_1 = \frac{1}{\mu}$, and $0 <\alpha \leq \min\Bigg\{\frac{1}{114}, \frac{\lp 1 -  \lamW^2 \rp^2}{432 \lamW^2}, \lp \frac{\lp 1 -  \lamW^2 \rp}{24\lamW} \rp^{2/3}, \frac{\mu}{6 L} , \frac{\mu^2 \lp 1 -  \lamW^2 \rp^2}{48 L^2 \lamW^2} \Bigg\}  $, and if $\mu \geq \frac{6 \sqrt{3} L}{n} \lp 1 + \frac{8 \lamW^2}{(1 - \lamW^2)}\rp$, then after further simplification we get $C_{\mu} \geq \frac{\mu}{4} > 0$, using which we get the desired result.
\end{IEEEproof}
Finally, we are ready to state the main theorem regarding the existence of $\epsilon$-KKT point. Specifically, we will bound $\frac{1}{n} \sumin \frac{1}{T} \sum_{t=1}^{T} \EE \lb \norm{\nabla u(\xh_i^t) +\hat{\w}_i^t }^2 \rb $ \eqref{opt}.
\begin{rem}
    It is remarked that combining the results of Lemma \ref{lem:sum_cnsnss_cmutatv} and \ref{lem:cmu} proves that the consensus is achieved 
\end{rem}
\begin{theorem} \label{the1}
Under the considered Assumptions \ref{assm:below}-\ref{assm:condn_on_graph}, if $ 0<\beta=\alpha^2<1, \mu \geq \frac{6 \sqrt{3} L}{n} \lp 1 + \frac{8 \lamW^2}{(1 - \lamW^2)}\rp$, and  $0 <\alpha \leq \min\Bigg\{\frac{1}{116}, \frac{\lp 1 -  \lamW^2 \rp^2}{432 \lamW^2}, \lp \frac{\lp 1 -  \lamW^2 \rp}{24\lamW} \rp^{2/3}, \frac{\mu}{6 L} , \frac{\mu^2 \lp 1 -  \lamW^2 \rp^2}{48 L^2 \lamW^2} \Bigg\}  $, then the mean squared stationary gap of the proposed D-MSSCA algorithm is upper bounded for all $T\geq 2$ as below:	
	\begin{align}
		\frac{1}{n} \sumin \frac{1}{T}\sum_{t= 1}^{T} \EE \lb \norm{\nabla u(\xh_i^t) + \hat{\w}_i^t}^2 \rb &\leq \lp \frac{8 L^2}{T}\lp 2 + \frac{9}{n}  + \frac{72 \lamW^2 }{n(1- \lamW^2)^2}   \rp \rp \frac{4}{\alpha \mu}\lp U(\bx^1) - U^\star \rp +\frac{48  \norm{\nabla u(\x_1^1)}^2 }{nT (1- \lamW^2)} P \nn \\
        &+ \frac{48 \bar{\sigma}^2}{nT b_0^2 (1- \lamW^2)} P + \frac{12\bar{\sigma}^2}{T n^2 b_0 \alpha^2} P + \frac{24\alpha^2 \bar{\sigma}^2 }{n} P \lp \frac{4\lamW^2}{T b_0 (1-\lamW^2)^2} + \frac{1}{n}\rp      \nn \\
		&  + \frac{72 \alpha^4\lamW^2 \bar{\sigma}^2 }{n(1 - \lamW^2)}  P + \frac{192\alpha^6\lamW^2 \bar{\sigma}^2 }{n(1 - \lamW^2)^2}  P, \label{eq:theo:1}
    \end{align}
	where, $P = \frac{8L^2 }{\mu^2} + \frac{36 L^2 }{n \mu^2} + \frac{288 L^2 \lamW^2}{n \mu^2 \lp 1 -  \lamW^2 \rp^2}+ 1 $.
\end{theorem}
\begin{IEEEproof}
	We will start the proof by using the optimality condition of \eqref{eq:alg_up_eq_1} to bound the mean squared stationary gap \eqref{opt}. 
 
    From the update equation\eqref{eq:alg_up_eq_1} and the definition of the surrogate function $\tf$, there exists $\hat{\w}_i^t \in \partial (h (\xh_i^t) + \ind_{\mathcal{X}} )$ for all $t \geq 1$ such that;    
     \begin{align}
         \nabla \fh_i \lp \xh_i^t, \x_i^t, \xi_i^t \rp  + (1 - \beta)\lp \v_i^{t-1} - \nfi(\x_i^{t-1},\xi_i^{t})\rp   + \pi_i^t + \hat{\w}_i^t  = 0. \nn
     \end{align}
     Further adding and subtracting $\nabla \fh_i \lp \x_i^t, \x_i^t, \xi_i^t \rp$, applying the definition of $\pi_i^t$ and update $\z_i^t$ \eqref{eq:alg_up_eq_3}, we obtain
    \begin{align}
         \nabla \fh_i \lp \xh_i^t, \x_i^t, \xi_i^t \rp - \nabla \fh_i \lp \x_i^t, \x_i^t, \xi_i^t \rp + \y_i^t  + \hat{\w}_i^t =0. \nn
    \end{align}
    By substituting $\hat{\w}_i^t = -\nabla \fh_i \lp \xh_i^t, \x_i^t, \xi_i^t \rp + \nabla \fh_i \lp \x_i^t, \x_i^t, \xi_i^t \rp - \y_i^t$, into the definition of mean squared stationary gap, we get
    \begin{align}
        \frac{1}{n} \sumin \frac{1}{T} \sum_{t=1}^{T} \EE \lb \norm{\nabla u(\xh_i^t) +\hat{\w}_i^t }^2 \rb =  \frac{1}{n} \sumin \frac{1}{T} \sum_{t=1}^{T} \EE \lb \norm{\nabla u(\xh_i^t) -\nabla \fh_i \lp \xh_i^t, \x_i^t, \xi_i^t \rp + \nabla \fh_i \lp \x_i^t, \x_i^t, \xi_i^t \rp - \y_i^t }^2 \rb. \nn
    \end{align}
    Now, by adding and subtracting $\nabla u(\x_i^t) - \nabla u(\bx^t)$ and separating the terms using the properties of the norm, we obtain:
    \begin{align}
        \frac{1}{n} &\sumin \frac{1}{T} \sum_{t=1}^{T} \EE \lb \norm{\nabla u(\xh_i^t) +\hat{\w}_i^t }^2 \rb \nn \\
        &\leq 
        \frac{4}{n} \sumin \frac{1}{T} \sumtT \Bigg( \EE \lb  \norm{\nabla u(\xh_i^t)  - \nabla u( \x_i^t)}^2 \rb + \EE \lb \norm{\nabla u(\x_i^t) - \nabla u( \bx^t)}^2 \rb + \EE \lb \norm{\nabla u( \bx^t) - \y_i^t}^2\rb^2  \nn \\
        &\quad+ \EE \lb \norm{\nabla \fh_i \lp \x_i^t, \x_i^t, \xi_i^t \rp  -\nabla \fh_i \lp \xh_i^t, \x_i^t, \xi_i^t \rp  }^2  \rb \Bigg), \nn \\
        &\leqtext{(i)} \frac{4}{n} \sumin \frac{1}{T}\sum_{t= 1}^{T} \Bigg( 2 L^2 \EE  \lb  \norm{\xh_i^t - \x_i^t}^2 \rb + L^2 \EE  \lb \norm{ \x_i^t - \bx^t}^2 \rb \Bigg) \nn \\
        &\qquad + \frac{4}{nT}\sum_{t= 1}^{T} \EE \lb \norm{(\bo\otimes \I_d) \lp \nabla u( \bx^t) -\bnu(\x^t) + \bnu(\x^t)- \by^t +\by^t \rp - \y^t}^2 \rb,  \nn \\
        &\leq \frac{8L^2}{nT} \sum_{t= 1}^{T} \EE  \lb  \norm{\xh^t - \x^t}^2 \rb + \frac{4L^2}{nT} \sumtT \EE  \lb \norm{ \x^t - (\bo\otimes \I_d) \bx^t}^2 \rb + \frac{4}{nT}\sum_{t= 1}^{T} \Bigg(3n \EE \lb \norm{\bnu(\x^t)- \by^t}^2 \rb \nn \\
        &\qquad \qquad+ 3\EE \lb \norm{(\bo\otimes \I_d) \lp \nabla u( \bx^t) -\bnu(\x^t) \rp}^2 \rb + 3 \EE \lb \norm{\y^t-(\bo\otimes \I_d) \by^t }^2 \rb \Bigg), \nn \\
        &\leqtext{(ii)} \frac{8L^2}{nT} \sum_{t= 1}^{T} \EE  \lb  \norm{\delb^t}^2 \rb + \frac{4L^2}{nT} \sumtT \EE  \lb \lp \theta^t\rp^2 \rb + \frac{12 L^2}{nT}\sum_{t= 1}^{T} \EE \lb \lp \theta^t \rp \rb + \frac{12}{T} \sum_{t= 1}^{T} \EE  \lb \norm{\bnu(\x^t)- \bz^t}^2 \rb \nn \\
        &\qquad + \frac{12}{nT} \sumtT \varepsilon^t. \nn
    \end{align}
    In (i), we have applied assumption~\ref{assm:smoothness}, and in (ii), we have applied \eqref{eq:sub:lem:basic_decentralized_result_2} and  \eqref{eq:sub:lem:basic_decentralized_result_3}. Now, by substituting the value of  $\sum_{t= 1}^{T} \EE  \lb \norm{\bnu(\x^t)- \bz^t}^2 \rb$, from Lemma~\ref{lem:HSGD}, we get
    \begin{align}
        \frac{1}{n} \sumin \frac{1}{T} \sum_{t=1}^{T} \EE \lb \norm{\nabla u(\xh_i^t) +\hat{\w}_i^t }^2 \rb &\leq \frac{8L^2}{nT} \sum_{t= 1}^{T} \EE  \lb  \norm{\delb^t}^2 \rb + \frac{16L^2}{nT} \sumtT \EE  \lb \lp \theta^t\rp^2 \rb + \frac{12}{nT} \sumtT \varepsilon^t \nn \\
        &+ \frac{12}{T} \lp \frac{\bar{\sigma}^2}{n^2 b_0 \beta} + \frac{2 \beta  \bar{\sigma}^2 T}{n^2} +  \frac{6L^2\alpha^2}{n^2\beta}  \sum_{t=1}^{T-1}  \EE \lb \norm{\delb^t}^2 \rb 
            + \frac{12 L^2}{\beta n^2} \sum_{t=1}^{T} \EE \lb \lp \theta^t \rp^2 \rb \rp. \nn
    \end{align}

Further combining the common terms and substituting the bound of $ \varepsilon^t $ from Lemma~\ref{lem:sumy}, we get
\begin{align}
    \frac{1}{n} &\sumin \frac{1}{T} \sum_{t=1}^{T} \EE \lb \norm{\nabla u(\xh_i^t) +\hat{\w}_i^t }^2 \rb \leq \frac{8L^2}{nT} \lp 1  +  \frac{9 \alpha^2}{n \beta} \rp \sum_{t= 1}^{T} \EE  \lb  \norm{\delb^t}^2 \rb + \frac{16L^2}{nT} \lp 1 + \frac{9 }{\beta n}\rp \sumtT \EE  \lb \lp \theta^t\rp^2 \rb  \nn \\
    &+ \frac{12}{T} \lp \frac{\bar{\sigma}^2}{n^2 b_0 \beta} + \frac{2 \beta  \bar{\sigma}^2 T}{n^2} \rp + \frac{12}{nT} \Bigg[  \frac{24 \alpha^2 \lamW^2 L^2}{\lp 1- \lamW^2 \rp^2}\lp 3 + 2 \beta \rp \sum_{t=1}^{T-1} \EE \lb \norm{\delb^t}^2  \rb +\frac{2}{1- \lamW^2}  \varepsilon^{1} \nn \\
    &+ \frac{24 \lamW^2 L^2}{\lp 1- \lamW^2 \rp^2} \lp 3 + 4 \beta  \rp \sum_{t=1}^{T} \EE \lb   \lp \theta^t \rp^2 \rb  +\frac{2 \lamW^2 \beta^2 \bar{\sigma}^2}{1- \lamW^2} \lp \frac{4}{\lp 1- \lamW^2 \rp b_0 \beta} + \frac{8 \beta T }{\lp 1- \lamW^2 \rp}  + 3 T \rp \Bigg]. \nn
\end{align}

Combining the common terms and substituting the bound of $\sumtT \EE \lb \lp \theta^t \rp^2\rb$ from Lemma~\ref{lem:sum_cnsnss_cmutatv}, we can further simplify as follows:
\begin{align}
    \frac{1}{n} &\sumin \frac{1}{T} \sum_{t=1}^{T} \EE \lb \norm{\nabla u(\xh_i^t) +\hat{\w}_i^t }^2 \rb \leq \frac{8L^2}{n}\lp  1  +  \frac{9 \alpha^2}{n \beta} + \frac{36 \alpha^2 \lamW^2 }{\lp 1- \lamW^2 \rp^2}\lp 3 + 2 \beta \rp  \rp \sum_{t= 1}^{T} \EE  \lb  \norm{\delb^t}^2 \rb \nn \\
    &+ \frac{12}{T} \lp \frac{\bar{\sigma}^2}{n^2 b_0 \beta} + \frac{2 \beta  \bar{\sigma}^2 T}{n^2} \rp + \frac{12}{nT} \Bigg[ \frac{2}{1- \lamW^2}  \varepsilon^{1} + \frac{2 \lamW^2 \beta^2 \bar{\sigma}^2}{1- \lamW^2} \lp \frac{4}{\lp 1- \lamW^2 \rp b_0 \beta} + \frac{8 \beta T }{\lp 1- \lamW^2 \rp}  + 3 T \rp \Bigg] \nn \\
    &+ \frac{16L^2}{nT} \lp  1 + \frac{9 }{\beta n} + \frac{18 \lamW^2}{\lp 1- \lamW^2 \rp^2} \lp 3 + 4 \beta  \rp \rp \frac{4 \alpha^2 \lamW^2 }{\lp 1 -  \lamW^2 \rp^2} \sum_{t=1}^{T-1} \EE \lb \norm{\delb^t}^2 \rb. \nn
\end{align}
Further using the bound of $\frac{1}{T}\sum_{t=1}^{T-1} \EE \lb \norm{\delb^t}^2 \rb $ obtained in Lemma~\ref{lem:cmu}, and applying \eqref{eq:lem:ref6}, we get
\begin{align}
    \frac{1}{n} &\sumin \frac{1}{T} \sum_{t=1}^{T} \EE \lb \norm{\nabla u(\xh_i^t) +\hat{\w}_i^t }^2 \rb \nn \\
    &\leq \frac{8L^2}{n} \lp 1 +\frac{9 \alpha^2}{n\beta}   +\frac{4 \alpha^2 \lamW^2 }{\lp 1 -  \lamW^2 \rp^2}\lp29+ 18\beta  +\frac{18}{n\beta} + \frac{36\lamW^2}{(1 - \lamW^2)^2 } \lp3 + 4\beta \rp   \rp  \rp \Bigg[ \frac{4n}{\alpha T \mu}U(\bx^1) \nn \\
    &- \frac{4n}{\alpha T \mu}U^\star  + \frac{24}{T \mu^2 \lp 1- \lamW^2 \rp} \norm{\nabla u(\x^1)}^2  + \frac{12 \bar{\sigma}^2}{2 T \mu^2} \Bigg(  \frac{1}{b_0 \alpha^2 n} + \frac{2 \alpha^2 T}{n} + \frac{4}{b_0^2 \lp 1- \lamW^2 \rp} \nn \\
    &+ \frac{2 \lamW^2 \alpha^4}{1- \lamW^2} \lp \frac{4}{\lp 1- \lamW^2 \rp b_0 \alpha^2} + \frac{8 \alpha^2 T }{\lp 1- \lamW^2 \rp}  + 3 T \rp \Bigg)  \Bigg] + \frac{12}{T} \lp \frac{\bar{\sigma}^2}{n^2 b_0 \beta} + \frac{2 \beta  \bar{\sigma}^2 T}{n^2} \rp \nn \\
    &+ \frac{12}{nT} \Bigg[ \frac{2 \lamW^2 \beta^2 \bar{\sigma}^2}{1- \lamW^2} \lp \frac{4}{\lp 1- \lamW^2 \rp b_0 \beta} + \frac{8 \beta T }{\lp 1- \lamW^2 \rp}  + 3 T \rp + \frac{2}{1- \lamW^2}  \lp \frac{2 \bar{\sigma}^2}{b_0^2} + 2 \norm{\nabla u(\x^1)}^2 \rp \Bigg]. \nn 
\end{align}
Further substituting $\beta = \alpha ^2$, considering $\alpha \leq \min \Bigg\{ \frac{1}{116},\frac{\lp 1 -  \lamW^2 \rp^2}{432\lamW^2}, \frac{\lp 1 -  \lamW^2 \rp^{2/3}}{8\lamW^{2/3}}, \frac{\lp 1 -  \lamW^2 \rp^2}{4 \lamW^2} \Bigg\}$, and rearranging, we get the desired result.
\end{IEEEproof}
\begin{rem}
    It can be noted from the Theorem~\ref{the1} that the mean squared stationary gap of the proposed D-MSSCA algorithm reaches to a steady state-error at a sublinear rate. Where the steady-state error is defined as $$\underset{T \rightarrow \infty}{\limsup} \frac{1}{n} \sumin \EE \lb \norm{\nabla u(\xh_i^t) +\hat{\w}_i^t }^2 \rb \leq  \frac{24\alpha^2 \bar{\sigma}^2 }{n^2} P $$
    From this expression, it is evident that the steady-state error can be reduced by selecting smaller values of $\alpha$ and $\beta$. Additionally, the error decreases as the number of nodes increases, which is expected since $n$ nodes function as $n$ oracles (SFO) for estimating the gradient.
\end{rem}
Finally, the next corollary provides the convergence rate in terms of the SFO complexity of the proposed D-MSSCA algorithm for fixed values of $\alpha, \beta$, and $b_0$.

\begin{cor} \label{col1} Under the conditions such that Theorem~\ref{the1} holds, if we further consider $\alpha=\mathcal{O}(T^{-1/3}), \beta=\mathcal{O}(T^{-2/3})$, and $b_0=\mathcal{O}(T^{1/3})$, The proposed D-MSSCA algorithm achieves an $\epsilon-$ KKT point in $\mathcal{O}(\epsilon^{-3/2})$ oracle calls. 
\end{cor} 
\begin{IEEEproof}
    By substituting $\alpha = T^{-1/3}, b_0 = T^{1/3}$ in \eqref{eq:theo:1}, we get
    \begin{align}
		\frac{1}{n} \sumin \frac{1}{T}\sum_{t= 1}^{T} \EE \lb \norm{\nabla u(\xh_i^t) + \hat{\w}_i^t}^2 \rb &\leq \lp \frac{8 L^2}{T^{2/3}}\lp 2 + \frac{9}{n}  + \frac{72 \lamW^2 }{n(1- \lamW^2)^2}   \rp \rp \frac{4}{ \mu}\lp U(\bx^1) - U^\star \rp \nn \\
        &+\frac{48  \norm{\nabla u(\x_1^1)}^2 }{nT (1- \lamW^2)} P + \frac{48 \bar{\sigma}^2}{n T^{5/3} (1- \lamW^2)} P + \frac{12\bar{\sigma}^2}{n^2 T^{2/3}} P \nn \\
		& + \frac{24 T^{-2/3} \bar{\sigma}^2 }{n} P \lp \frac{4\lamW^2}{ T^{4/3} (1-\lamW^2)^2} + \frac{1}{n}\rp  + \frac{72 T^{-4/3} \lamW^2 \bar{\sigma}^2 }{n(1 - \lamW^2)}  P \nn \\
        &+ \frac{192 T^{-6/3} \lamW^2 \bar{\sigma}^2 }{n(1 - \lamW^2)^2}  P, \nn
    \end{align}
    In order to reach $\epsilon-$KKT point, we require
    \begin{align}
        \frac{1}{T^{2/3}}&\Bigg[\frac{32 L^2}{\mu} \lp 2 + \frac{9}{n}  + \frac{72 \lamW^2 }{n(1- \lamW^2)^2}   \rp \lp U(\bx^1) - U^\star \rp +\frac{48 P \norm{\nabla u(\x_1^1)}^2 }{nT^{1/3} (1- \lamW^2)}+ \frac{48 P \bar{\sigma}^2}{n T (1- \lamW^2)}  + \frac{12P\bar{\sigma}^2}{n^2} \nn \\
		& + \frac{24 P \bar{\sigma}^2 }{n} \lp \frac{4\lamW^2}{ T^{4/3} (1-\lamW^2)^2} + \frac{1}{n}\rp  + \frac{72 P\lamW^2 \bar{\sigma}^2 }{n T^{2/3}(1 - \lamW^2)} + \frac{192 T^{-6/3} \lamW^2 \bar{\sigma}^2 }{n(1 - \lamW^2)^2}  P \Bigg] < \epsilon,
    \end{align}
    which gives $T=\mathcal{O}(\epsilon^{-3/2})$.
\end{IEEEproof}

The SFO complexity of Corollary \ref{col1} matches that of DEEPSTORM \cite{mancino2023proximal} and ProxGT-SR-O/E \cite{xin2021stochastic}. It should be noted that, unlike ProxGT-SR-O/E, which requires large batch sizes, our algorithm is batchless and uses one sample at each iteration. Also, this rate matches the SFO complexity lower bound for centralized unconstrained stochastic non-convex optimization problems. Also, this rate matches the SFO complexity lower bound for centralized unconstrained stochastic non-convex optimization problems.

\section{Experimental data and results} \label{resultss}
\begin{figure}
      \centering
     \includegraphics[width = 0.5\textwidth]{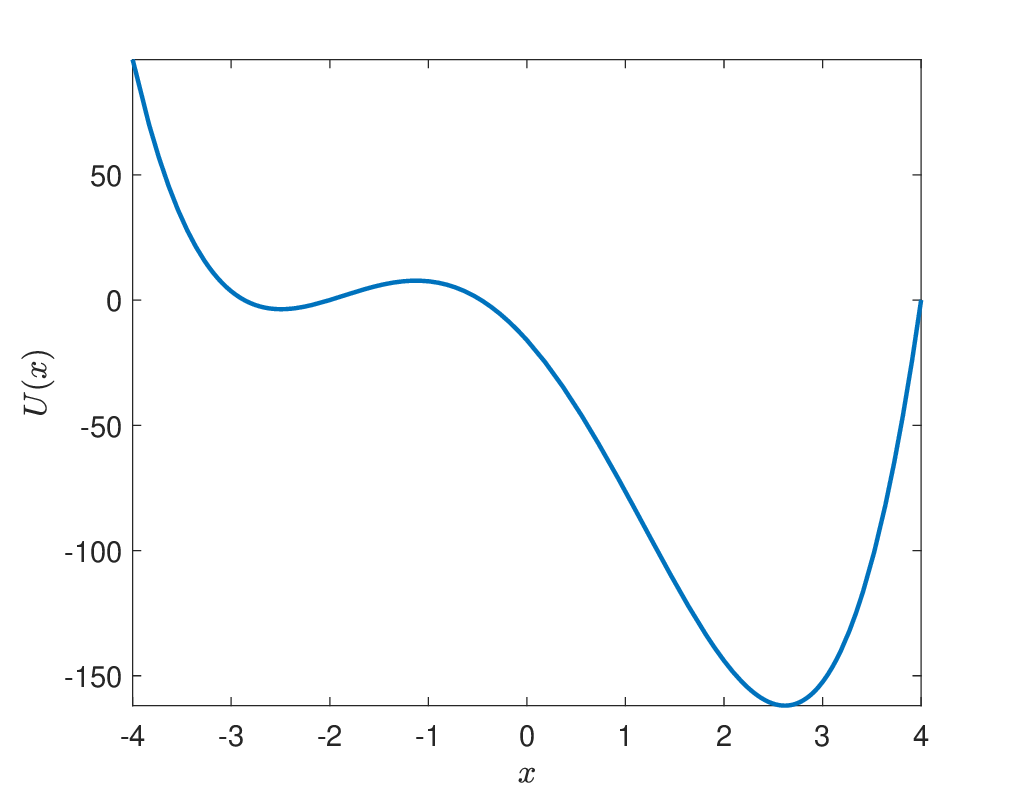}
    \caption[center]{Golbal objective function $U(x)$.}
                   \label{fig4}
    \end{figure}

\begin{figure}
      \centering
     \includegraphics[width = 0.5\textwidth]{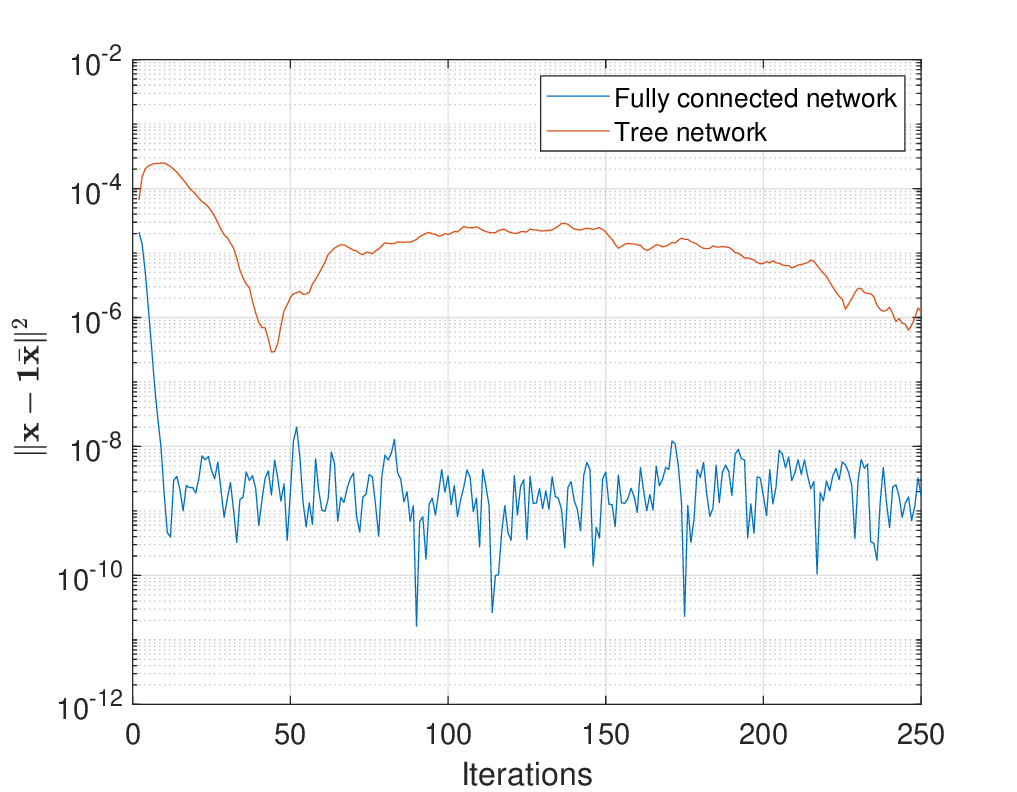}
    \caption[center]{Evolution of residual $\norm{\x - \bo \bx }^2$ over different networks.}
                   \label{fig3}
    \end{figure}

\begin{figure}
      \centering
     \includegraphics[width = 0.5\textwidth]{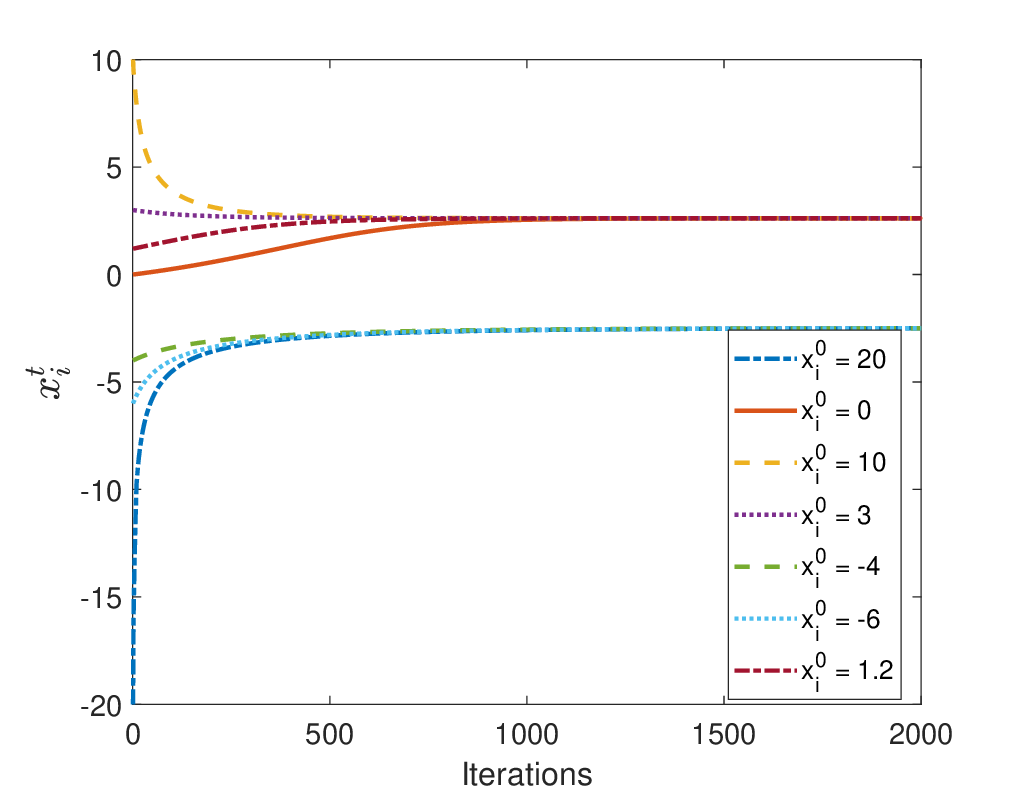}
    \caption[center]{Evolution of local variable for varying initialization.}
                   \label{fig1}
    \end{figure}

\begin{figure}
      \centering
     \includegraphics[width = 0.5\textwidth]{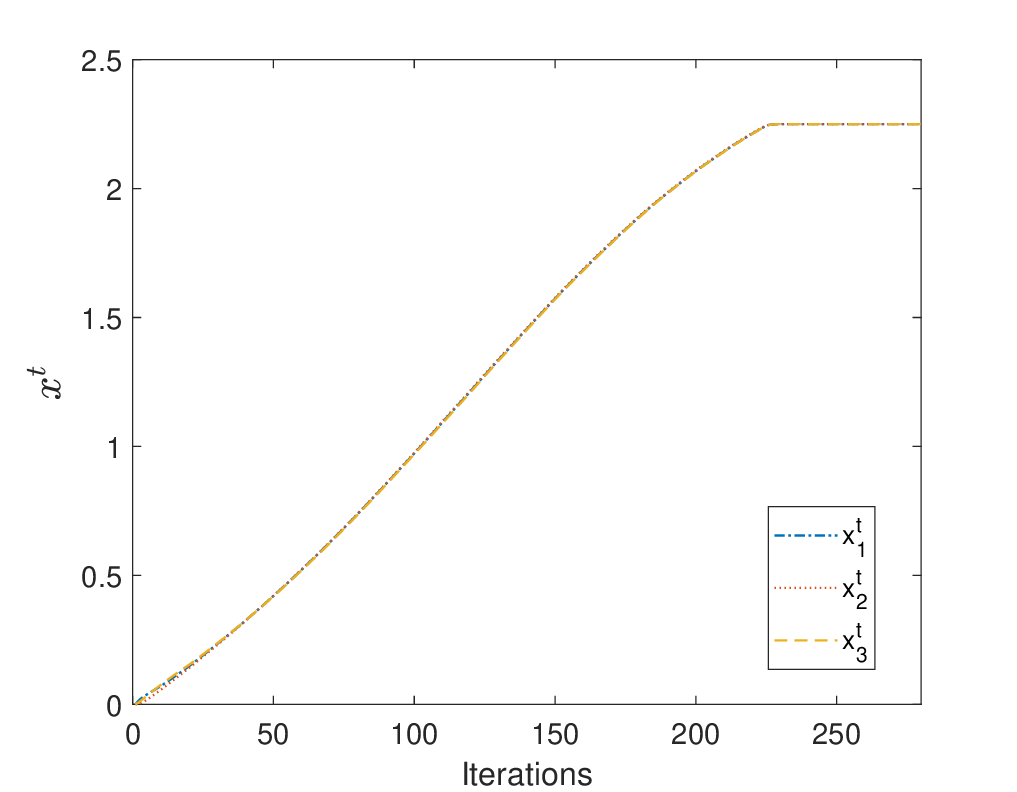}
    \caption[center]{Evolution of local variable when each node is initialized at $\x_i^t = 0$ given a global constraint $\abs{x_i} \leq 2.25$.}
                   \label{fig2}
    \end{figure}

In this section, we will demonstrate the applicability of D-MSSCA. Let us consider a simple distributed optimization problem, which is a stochastic version of the synthetic problem in \cite{tatarenko2017non, zheng2023distributed} over a network of $n = 3$ nodes:
\begin{align}
    U(x) &= \min \sum_{i = 1}^3 \EE\lb f_i(x,\xi_i)\rb 
\end{align}
Each local objective function $f_i$ is defined as
\begin{align}
		f_1(x,\xi_1) = \begin{cases}
			(x^3 - 16x )(x+2) + n_1x, &  |x|\leq 10  \\
			4248x - 32400 + n_1x, &  x > 10 \\
            -3112x - 25040+ n_1x, & x < -10
		\end{cases}      \\
        f_2(x,\xi_2) = \begin{cases}
			(0.5x^3 + x^2 )(x-4)+ n_2x, &  |x|\leq 10  \\
			1620x - 12600+ n_2x, &  x > 10 \\
            -2220x - 16600+ n_2x, & x < -10
		\end{cases}   \\
        f_3(x,\xi_3) = \begin{cases}
			(x^3 -16x )(x+2)+ n_3x, &  |x|\leq 10  \\
			288x - 2016+ n_3x, &  x > 10 \\
            228x - 2624+ n_3x, & x < -10
		\end{cases} 
 \end{align}
where $\xi_i = n_i \sim \mathcal{N}(0, 1)$. The objective $U(x)$ in the interval $[-4,4]$ is shown in Fig. \ref{fig4}. As the objective function is non-convex, and its surrogate can be constructed using \eqref{def:tf}. We use D-MSSCA to first demonstrate the effect of communication topology in Fig \ref{fig3}. We observe that the fully connected network performs better than the Tree network, which is how the behavior is expected. Next,  setting $\alpha =0.8 $ and $\beta = 0.16$ and $\mu = 5000$,  we plot the evolution of local variable $x_i^t$ with different initial values in Fig \ref{fig1} for a fully connected graph with $\lamW = 0.5$. We observe the nodes converge to local minima. Finally, we plot the evolution of local variables of each node given a global constraint set $\abs{x_i} \leq 2.25$, when all the nodes are initialized at $x_i^t = 0$ in Fig \ref{fig2}. It can be observed that all the local variables get as close as possible to the true minima.

\section{Conclusion and Future work}

In this work, we consider decentralized consensus stochastic non-convex optimization to minimize the sum of non-convex (possibly smooth) and convex (possibly non-smooth) cost functions over a network of nodes. While this problem is well studied, comprehensive convergence analysis under the stochastic SCA rubric has remained an open problem. We have proposed and analyzed D-MSSCA that achieves the optimal rate of $\mathcal{O}(\epsilon^{-3/2})$. The rate matches the SFO rate of the state-of-the-art decentralized gradient-based algorithm while processing a single sample at each iteration. The algorithm uses strong convex surrogates and leverages recursive momentum-based updates at each node, achieving faster convergence. The applicability of D-MSSCA is demonstrated on a synthetic stochastic problem. One interesting future direction of this paper is under investigation, wherein we simplify the optimization problem \eqref{eq:alg_up_eq_1} in Algorithm \ref{alg:DMSCA} by linearizing the global constraint $g$ making it \eqref{eq:alg_up_eq_1} easier to solve than proximal-based methods. By reducing the complexity of each iteration, this modification could improve convergence rates and expand the applicability of the D-MSSCA algorithm to broader classes of optimization problems, especially those where proximal methods face scalability and efficiency challenges.

\bibliographystyle{IEEEtran} 
\bibliography{IEEEabrv,cittt}

\begin{thebibliography}{10}
\providecommand{\url}[1]{#1}
\csname url@samestyle\endcsname
\providecommand{\newblock}{\relax}
\providecommand{\bibinfo}[2]{#2}
\providecommand{\BIBentrySTDinterwordspacing}{\spaceskip=0pt\relax}
\providecommand{\BIBentryALTinterwordstretchfactor}{4}
\providecommand{\BIBentryALTinterwordspacing}{\spaceskip=\fontdimen2\font plus
\BIBentryALTinterwordstretchfactor\fontdimen3\font minus
  \fontdimen4\font\relax}
\providecommand{\BIBforeignlanguage}[2]{{%
\expandafter\ifx\csname l@#1\endcsname\relax
\typeout{** WARNING: IEEEtran.bst: No hyphenation pattern has been}%
\typeout{** loaded for the language `#1'. Using the pattern for}%
\typeout{** the default language instead.}%
\else
\language=\csname l@#1\endcsname
\fi
#2}}
\providecommand{\BIBdecl}{\relax}
\BIBdecl

\bibitem{bottou2018optimization}
L.~Bottou, F.~E. Curtis, and J.~Nocedal, ``Optimization methods for large-scale
  machine learning,'' \emph{SIAM review}, vol.~60, no.~2, pp. 223--311, 2018.

\bibitem{bianchi2012convergence}
P.~Bianchi and J.~Jakubowicz, ``Convergence of a multi-agent projected
  stochastic gradient algorithm for non-convex optimization,'' \emph{IEEE
  Trans. Autom. Control}, vol.~58, no.~2, pp. 391--405, 2012.

\bibitem{swenson2022distributed}
B.~Swenson, R.~Murray, H.~V. Poor, and S.~Kar, ``Distributed stochastic
  gradient descent: Nonconvexity, nonsmoothness, and convergence to local
  minima,'' \emph{Journal of Machine Learning Research}, vol.~23, no. 328, pp.
  1--62, 2022.

\bibitem{wang2021distributed}
Z.~Wang, J.~Zhang, T.-H. Chang, J.~Li, and Z.-Q. Luo, ``Distributed stochastic
  consensus optimization with momentum for nonconvex nonsmooth problems,''
  \emph{IEEE Trans. Signal Process}, vol.~69, pp. 4486--4501, 2021.

\bibitem{xin2021stochastic}
R.~Xin, S.~Das, U.~A. Khan, and S.~Kar, ``A stochastic proximal gradient
  framework for decentralized non-convex composite optimization:
  Topology-independent sample complexity and communication efficiency,''
  \emph{arXiv preprint arXiv:2110.01594}, 2021.

\bibitem{yan2023compressed}
Y.~Yan, J.~Chen, P.-Y. Chen, X.~Cui, S.~Lu, and Y.~Xu, ``Compressed
  decentralized proximal stochastic gradient method for nonconvex composite
  problems with heterogeneous data,'' in \emph{International Conference on
  Machine Learning}.\hskip 1em plus 0.5em minus 0.4em\relax PMLR, 2023, pp.
  39\,035--39\,061.

\bibitem{mancino2023proximal}
G.~Mancino-Ball, S.~Miao, Y.~Xu, and J.~Chen, ``Proximal stochastic recursive
  momentum methods for nonconvex composite decentralized optimization,'' in
  \emph{Proceedings of the AAAI Conference on Artificial Intelligence},
  vol.~37, no.~7, 2023, pp. 9055--9063.

\bibitem{zheng2023distributed}
L.~Zheng, H.~Li, J.~Li, Z.~Wang, Q.~L{\"u}, Y.~Shi, H.~Wang, T.~Dong, L.~Ji,
  and D.~Xia, ``A distributed nesterov-like gradient tracking algorithm for
  composite constrained optimization,'' \emph{IEEE Trans. Signal. Inf. Process
  Netw.}, vol.~9, pp. 60--73, 2023.

\bibitem{di2019distributed}
P.~Di~Lorenzo and S.~Scardapane, ``Distributed stochastic nonconvex
  optimization and learning based on successive convex approximation,'' in
  \emph{2019 53rd Asilomar Conference on Signals, Systems, and
  Computers}.\hskip 1em plus 0.5em minus 0.4em\relax IEEE, 2019, pp. 1--5.

\bibitem{scutari2013decomposition}
G.~Scutari, F.~Facchinei, P.~Song, D.~P. Palomar, and J.-S. Pang,
  ``Decomposition by partial linearization: Parallel optimization of
  multi-agent systems,'' \emph{IEEE Trans. Signal Process}, vol.~62, no.~3, pp.
  641--656, 2013.

\bibitem{scutari2016parallel}
G.~Scutari, F.~Facchinei, and L.~Lampariello, ``Parallel and distributed
  methods for constrained nonconvex optimization—part {I}: Theory,''
  \emph{IEEE Trans. Signal Process.}, vol.~65, no.~8, pp. 1929--1944, 2016.

\bibitem{scutari2016paralel}
G.~Scutari, F.~Facchinei, L.~Lampariello, S.~Sardellitti, and P.~Song,
  ``Parallel and distributed methods for constrained nonconvex
  optimization-part {II}: Applications in communications and machine
  learning,'' \emph{IEEE Trans. Signal Process.}, vol.~65, no.~8, pp.
  1945--1960, 2016.

\bibitem{yang2016parallel}
Y.~Yang, G.~Scutari, D.~P. Palomar, and M.~Pesavento, ``A parallel
  decomposition method for nonconvex stochastic multi-agent optimization
  problems,'' \emph{IEEE Trans. Signal Process.}, vol.~64, no.~11, pp.
  2949--2964, 2016.

\bibitem{liu2019stochastic}
A.~Liu, V.~K. Lau, and B.~Kananian, ``Stochastic successive convex
  approximation for non-convex constrained stochastic optimization,''
  \emph{IEEE Trans. Signal Process.}, vol.~67, no.~16, pp. 4189--4203, 2019.

\bibitem{liu2018online}
A.~Liu, V.~K. Lau, and M.-J. Zhao, ``Online successive convex approximation for
  two-stage stochastic nonconvex optimization,'' \emph{IEEE Trans. Signal
  Process.}, vol.~66, no.~22, pp. 5941--5955, 2018.

\bibitem{ye2019stochastic}
C.~Ye and Y.~Cui, ``Stochastic successive convex approximation for general
  stochastic optimization problems,'' \emph{IEEE Wireless Commun. Lett.},
  vol.~9, no.~6, pp. 755--759, 2019.

\bibitem{liu2019two}
A.~Liu, X.~Chen, W.~Yu, V.~K. Lau, and M.-J. Zhao, ``Two-timescale hybrid
  compression and forward for massive {MIMO} aided {C-RAN},'' \emph{IEEE Trans.
  Signal Process.}, vol.~67, no.~9, pp. 2484--2498, 2019.

\bibitem{liu2021two}
A.~Liu, R.~Yang, T.~Q. Quek, and M.-J. Zhao, ``Two-stage stochastic
  optimization via primal-dual decomposition and deep unrolling,'' \emph{IEEE
  Trans. Signal Process.}, vol.~69, pp. 3000--3015, 2021.

\bibitem{mokhtari2017large}
A.~Mokhtari, A.~Koppel, G.~Scutari, and A.~Ribeiro, ``Large-scale nonconvex
  stochastic optimization by doubly stochastic successive convex
  approximation,'' in \emph{IEEE ICASSP}, 2017, pp. 4701--4705.

\bibitem{koppel2018parallel}
A.~Koppel, A.~Mokhtari, and A.~Ribeiro, ``Parallel stochastic successive convex
  approximation method for large-scale dictionary learning,'' in \emph{IEEE
  ICASSP}, 2018, pp. 2771--2775.

\bibitem{mokhtari2020high}
A.~Mokhtari and A.~Koppel, ``High-dimensional nonconvex stochastic optimization
  by doubly stochastic successive convex approximation,'' \emph{IEEE Trans.
  Signal Process.}, vol.~68, pp. 6287--6302, 2020.

\bibitem{idrees2021practical}
B.~M. Idrees, J.~Akhtar, and K.~Rajawat, ``Practical precoding via asynchronous
  stochastic successive convex approximation,'' \emph{IEEE Trans. Signal
  Process}, vol.~69, pp. 4177--4191, 2021.

\bibitem{idrees2024constrained}
B.~M. Idrees, L.~Arora, and K.~Rajawat, ``Constrained stochastic recursive
  momentum successive convex approximation,'' \emph{arXiv preprint
  arXiv:2404.11790}, 2024.

\bibitem{di2016next}
P.~Di~Lorenzo and G.~Scutari, ``Next: In-network nonconvex optimization,''
  \emph{IEEE Trans. Signal. Inf. Process Netw.}, vol.~2, no.~2, pp. 120--136,
  2016.

\bibitem{xin2021hybrid}
R.~Xin, U.~Khan, and S.~Kar, ``A hybrid variance-reduced method for
  decentralized stochastic non-convex optimization,'' in \emph{International
  Conference on Machine Learning}.\hskip 1em plus 0.5em minus 0.4em\relax PMLR,
  2021, pp. 11\,459--11\,469.

\bibitem{cutkosky2019momentum}
A.~Cutkosky and F.~Orabona, ``Momentum-based variance reduction in non-convex
  sgd,'' in \emph{Advances in Neural Information Processing Systems}, 2019, pp.
  15\,236--15\,245.

\bibitem{tran2019hybrid}
Q.~Tran-Dinh, N.~H. Pham, D.~T. Phan, and L.~M. Nguyen, ``Hybrid stochastic
  gradient descent algorithms for stochastic nonconvex optimization,''
  \emph{arXiv preprint arXiv:1905.05920}, 2019.

\bibitem{nguyen2017sarah}
L.~M. Nguyen, J.~Liu, K.~Scheinberg, and M.~Tak{'a}{v{c}}, ``Sarah a novel
  method for machine learning problems using stochastic recursive gradient,''
  in \emph{International Conference on Machine Learning}.\hskip 1em plus 0.5em
  minus 0.4em\relax PMLR, 2017, pp. 2613--2621.

\bibitem{sharma2024optimized}
S.~D. Sharma and K.~Rajawat, ``Optimized gradient tracking for decentralized
  online learning,'' \emph{IEEE Trans. Signal Process}, vol.~72, pp.
  1443--1459, 2024.

\bibitem{qu2017harnessing}
G.~Qu and N.~Li, ``Harnessing smoothness to accelerate distributed
  optimization,'' \emph{IEEE Trans. Control Netw. Syst.}, vol.~5, no.~3, pp.
  1245--1260, 2017.

\bibitem{xin2021improved}
R.~Xin, U.~A. Khan, and S.~Kar, ``An improved convergence analysis for
  decentralized online stochastic non-convex optimization,'' \emph{IEEE Trans.
  Signal Process}, vol.~69, pp. 1842--1858, 2021.

\bibitem{xin2022fast}
------, ``Fast decentralized nonconvex finite-sum optimization with recursive
  variance reduction,'' \emph{SIAM Journal on Optimization}, vol.~32, no.~1,
  pp. 1--28, 2022.

\bibitem{tatarenko2017non}
T.~Tatarenko and B.~Touri, ``Non-convex distributed optimization,'' \emph{IEEE
  Trans. Autom. Control}, vol.~62, no.~8, pp. 3744--3757, 2017.

\end{thebibliography}

\pagebreak

\appendices
\section{Proof of Lemma~\ref{lem:consnss_bound}}\label{pr:lem:consnss_bound}
\begin{IEEEproof}
    From the definition of $\theta^t$\eqref{eq:def_consensus_error} and $\x^t$-update \eqref{eq:alg_up_cmpct_eq_2}, we have
    \begin{align}
        \lp \theta^t \rp^2 &= \norm{\x^t - \frac{1}{n}\lp \bo_n \otimes \I_d \rp \bx^t}^2 \nn,\\
        &=\norm{\lp \I- \frac{1}{n} \bo_n \bo_n^\mathsf{T}\otimes \I_d \rp \Wu \lp \x^{t-1} + \alpha \lp \xh^{t-1} - \x^{t-1} \rp \rp}^2 \nn, \\
        &\eqtext{(a)}\norm{\lp \Wu- \frac{1}{n} \bo_n \bo_n^\mathsf{T}\otimes \I_d \rp \x^{t-1} + \alpha \lp \Wu- \frac{1}{n} \bo_n \bo_n^\mathsf{T}\otimes \I_d \rp \lp \xh^{t-1} - \x^{t-1} \rp}^2 \nn, \\
        &\leqtext{(b)} \lp 1 + \eta_1 \rp \norm{\Wu \x^{t-1}- \lp \frac{1}{n} \bo_n \bo_n^\mathsf{T}\otimes \I_d \rp \x^{t-1}}^2 \nn \\
        &\quad + \lp 1 + \frac{1}{\eta_1} \rp \alpha^2 \lp \lambda_{\max} \lp \Wu- \frac{1}{n} \bo_n \bo_n^\mathsf{T}\otimes \I_d\rp \rp^2 \norm{ \xh^{t-1} - \x^{t-1} }^2 \nn, \\
        &\leqtext{(c)} \lp 1 + \eta_1 \rp \lamW^2 \norm{\x^{t-1}- \lp \frac{1}{n} \bo_n \bo_n^\mathsf{T}\otimes \I_d \rp \x^{t-1}}^2 + \lp 1 + \frac{1}{\eta_1} \rp \alpha^2 \lamW^2 \norm{ \xh^{t-1} - \x^{t-1} }^2 \nn,\\
        &= \lp 1 + \eta_1 \rp \lamW^2 \lp \theta^{t-1}\rp^2 + \lp 1 + \frac{1}{\eta_1} \rp \alpha^2 \lamW^2 \norm{\delb^{t-1} }^2 \nn.
    \end{align}
     In (a), we have applied the property of $\W$, i.e., $\bo^\mathsf{T} \W
     = \bo^\mathsf{T}$. In (b), Young's inequality and the property of norm are used, i.e., for any square matrix $\A$ and vector $\x$ of compatible size, $\norm{\A\x}^2 \leq \lambda_{\max}(\A)^2 \norm{\x}^2$. In (c), we have applied \eqref{eq:sub:lem:basic_decentralized_result_1} and the definition of $\lamW$. 
\end{IEEEproof}

\section{Proof of Lemma~\ref{lem:bndnss_y}} \label{pr:lem:bndnss_y}
\begin{IEEEproof} From the definition of $\varepsilon^t$ and \eqref{eq:alg_up_cmpct_eq_4} we have,
    \begin{align}
         \varepsilon^t &= \EE\lb \norm{\Wu \lp \y^{t-1} + \z^{t} - \z^{t-1} \rp - \frac{1}{n}\lp \bo_n \bo_n^\mathsf{T} \otimes \I_d \rp \Wu \lp \y^{t-1} + \z^{t} - \z^{t-1} \rp}^2\rb \nn, \\
         &\eqtext{(i)}\EE\lb \norm{\Wu \lp \y^{t-1} + \z^{t} - \z^{t-1} \rp - \frac{1}{n}\lp \bo_n \bo_n^\mathsf{T} \otimes \I_d \rp \lp \y^{t-1} + \z^{t} - \z^{t-1} \rp}^2\rb \nn,\\
         &=\EE \Bigg[ \norm{\Wu \y^{t-1} - \frac{1}{n}\lp \bo_n \bo_n^\mathsf{T} \otimes \I_d \rp \y^{t-1} }^2 + \norm{\lp \Wu - \frac{1}{n} \bo_n\bo_n^\mathsf{T} \otimes \I_d \rp \lp \z^t - \z^{t-1} \rp }^2\nn \\
         &\quad + \left\langle \Wu \y^{t-1} - \frac{1}{n}\lp \bo_n\bo_n^\mathsf{T} \otimes \I_d \rp \y^{t-1}, \lp \Wu - \frac{1}{n} \bo_n\bo_n^\mathsf{T} \otimes \I_d \rp \lp \z^t - \z^{t-1} \rp \right \rangle \Bigg] \nn, \\
         &\leqtext{\eqref{eq:sub:lem:basic_decentralized_result_1}}
          \lamW^2 \EE \Bigg[\norm{\y^{t-1} - \frac{1}{n}\lp \bo_n \otimes \I_d \rp \by^{t-1} }^2 \Bigg]+ \lamW^2 \EE \Bigg[\norm{\z^t - \z^{t-1}  }^2 \Bigg]\nn \\
         &\quad + 2\EE \Bigg[\left\langle \Wu \y^{t-1} - \frac{1}{n}\lp \bo_n\bo_n^\mathsf{T} \otimes \I_d \rp \y^{t-1}, \lp \Wu - \frac{1}{n} \bo_n \bo_n^\mathsf{T}\otimes \I_d \rp \lp \z^t - \z^{t-1} \rp \right \rangle \Bigg]. \label{eq:lem:pr_1}
    \end{align}
    Where (i), uses the doubly stochastic matrix property of $\W$. Next, each term in the above equation will be simplified separately, starting with $\norm{\z^t - \z^{t-1}}^2$ we have,
    \begin{align}
       \norm{\z^t - \z^{t-1}}^2 &= \sumin \norm{\z_i^t - \z_i^{t-1}}^2 \nn, \\
       &\eqtext{\eqref{eq:alg_up_eq_3}} \sumin \norm{\nfi(\x_i^{t},\xi_i^{t}) + (1-\beta) \lp \z_i^{t-1} - \nfi(\x_i^{t-1},\xi_i^{t})\rp - \z_i^{t-1}}^2 \nn,\\
       &= \sumin \norm{\nfi(\x_i^t, \xi_i^t) - \nfi (\x_i^{t-1},\xi_i^t) - \beta\z_i^{t-1}  + \beta \nfi (\x_i^{t-1}, \xi_i^t)}^2 \label{eq:lem:pr_2}.
    \end{align}
  After adding and subtracting $\beta \nuu_i(\x^{t-1})$, we observe that
    \begin{align}
        &\norm{\z^t - \z^{t-1}}^2 \nn \\
        &~= \sumin \norm{\nfi(\x_i^t, \xi_i^t) - \nfi (\x_i^{t-1},\xi_i^t) - \beta \lp \z_i^{t-1} - \nuu_i(\x^{t-1}) \rp  + \beta \nfi (\x_i^{t-1}, \xi_i^t) -  \beta \nuu_i(\x^{t-1}) }^2 \nn, \\
        &\quad \leq \sumin 3\norm{\nfi(\x_i^t, \xi_i^t) - \nfi (\x_i^{t-1},\xi_i^t)}^2  + \sumin 3 \beta^2 \norm{\lp \z_i^{t-1} - \nuu_i(\x^{t-1}) \rp}^2  \nn \\
        &\qquad + \sumin 3 \beta^2 \norm{ \nfi (\x_i^{t-1}, \xi_i^t) -  \nuu_i(\x^{t-1}) }^2 \nn.
    \end{align}
    Taking expectation on both sides and further using Assumptions~\ref{assm:bounded_var} and \ref{assm:smoothness} we get
    \begin{align}
        \EE \lb \norm{\z^t - \z^{t-1}}^2 \rb \leq 3 L^2  \EE \lb \norm{\x^t -\x^{t-1}}^2 \rb  + 3 \beta^2  \EE \lb \norm{\z^{t-1} - \nuu(\x^{t-1})}^2 \rb + 3 \beta^2 \bar{\sigma}^2 \label{eq:lem:pr_3}.
    \end{align}
    Considering conditional expectation, we can further simplify the last term of \eqref{eq:lem:pr_1} as 
    \begin{align}
        2\EE &\Bigg[\left\langle \Wu \y^{t-1} - \frac{1}{n}\lp \bo_n\bo_n^\mathsf{T} \otimes \I_d \rp \y^{t-1}, \lp \Wu - \frac{1}{n} \bo_n \bo_n^\mathsf{T}\otimes \I_d \rp \lp \z^t - \z^{t-1} \rp \right \rangle \Bigg] \nn, \\
        &= 2\EE \Bigg[ \EE \lb \left\langle \Wu \y^{t-1} - \frac{1}{n}\lp \bo_n\bo_n^\mathsf{T} \otimes \I_d \rp \y^{t-1}, \lp \Wu - \frac{1}{n} \bo_n \bo_n^\mathsf{T}\otimes \I_d \rp \lp \z^t - \z^{t-1} \rp \right \rangle \Bigg|  \cH_k \rb \Bigg] \nn, \\
        &= 2\EE \Bigg[ \left\langle \Wu \y^{t-1} - \frac{1}{n}\lp \bo_n\bo_n^\mathsf{T} \otimes \I_d \rp \y^{t-1}, \lp \Wu - \frac{1}{n} \bo_n \bo_n^\mathsf{T}\otimes \I_d \rp \EE \lb \z^t - \z^{t-1} |  \cH_k  \rb \right \rangle  \Bigg]. \nn
    \end{align}
    Moreover, using \eqref{eq:lem:pr_2} and Assumption~\ref{assm:exp_equal_true}, we have $\EE \lb \z^t - \z^{t-1} |  \cH_k  \rb = \nuu(\x^t)-\nuu(\x^t)-\beta(\z^{t-1}-\nuu(\x^{t-1}))$. Substituting which along with Cauchy Schwartz inequality, gives 
    \begin{align}
        2\EE &\Bigg[\left\langle \Wu \y^{t-1} - \frac{1}{n}\lp \bo_n\bo_n^\mathsf{T} \otimes \I_d \rp \y^{t-1}, \lp \Wu - \frac{1}{n} \bo_n \bo_n^\mathsf{T}\otimes \I_d \rp \lp \z^t - \z^{t-1} \rp \right \rangle \Bigg] \nn \\
        &= 2\EE \Bigg[ \norm{ \Wu \y^{t-1} - \frac{1}{n}\lp \bo_n\bo_n^\mathsf{T} \otimes \I_d \rp \y^{t-1}} \times \nn \\
        &\quad \norm{\lp \Wu - \frac{1}{n} \bo_n \bo_n^\mathsf{T}\otimes \I_d \rp \lp \nuu(\x^t)-\nuu(\x^t)-\beta(\z^{t-1}-\nuu^{t-1}) \rp }  \Bigg] \nn, \\
        &\leqtext{\eqref{eq:sub:lem:basic_decentralized_result_1}} 2\EE \Bigg[\lamW^2 \norm{ \y^{t-1} - \frac{1}{n}\lp \bo_n \otimes \I_d \rp \by^{t-1}} \norm{\nuu(\x^t)-\nuu(\x^t)-\beta(\z^{t-1}-\nuu^{t-1})} \Bigg]. \nn
    \end{align}
   Young's inequality along with few mathematical simplifications gives,
   \begin{align}
        \EE &\Bigg[\left\langle \Wu \y^{t-1} - \frac{1}{n}\lp \bo_n\bo_n^\mathsf{T} \otimes \I_d \rp \y^{t-1}, \lp \Wu - \frac{1}{n} \bo_n \bo_n^\mathsf{T}\otimes \I_d \rp \lp \z^t - \z^{t-1} \rp \right \rangle \Bigg] \nn \\
        &\leq \lamW^2 \gamma_1 \EE \lb \norm{ \y^{t-1} - \frac{1}{n}\lp \bo_n \otimes \I_d \rp \by^{t-1}}^2 \rb + \frac{2 \lamW^2 L^2}{\gamma_1} \EE \lb \norm{\x^t-\x^{t-1}}^2 \rb \nn \\
        &+ \frac{2 \lamW^2 \beta^2}{\gamma_1} \EE \lb \norm{(\z^{t-1}-\nuu^{t-1})}^2 \rb. \label{eq:lem:pr_4}
    \end{align}
    By substituting equations \eqref{eq:lem:pr_3}, \eqref{eq:lem:pr_4} in \eqref{eq:lem:pr_1}, we can write
    \begin{align}
         \varepsilon^t &\leq 
          \lamW^2 \lp 1 + \gamma_1\rp \varepsilon^{t-1} + \lamW^2 \beta^2 \lp \frac{2 }{\gamma_1} + 3 \rp \upsilon^{t-1}   + 3 \lamW^2 \beta^2 \bar{\sigma}^2  \nn \\
          &\quad + \lamW^2 L^2 \lp 3 + \frac{2}{\gamma_1}\rp \EE \lb \norm{\x^t - \bo \otimes \I_d \bx^{t}+ \bo \otimes \I_d \bx^{t}- \bo \otimes \I_d \bx^{t-1}+\bo \otimes \I_d \bx^{t-1}-\x^{t-1}}^2 \rb.  \nn 
    \end{align}
    We can further simplify the last term as follows
    \begin{align}
        \norm{\x^t - \x^{t-1}}^2 &= \norm{\x^t - \bo \otimes \I_d \bx^{t}+ \bo \otimes \I_d \bx^{t}- \bo \otimes \I_d \bx^{t-1}+\bo \otimes \I_d \bx^{t-1}-\x^{t-1}}^2   \nn, \\
        &\leq  3 \norm{\x^t - \bo \otimes \I_d \bx^{t}}^2 + 3\norm{\bo \otimes \I_d \bx^{t}- \bo \otimes \I_d \bx^{t-1}}^2 + 3\norm{\bo \otimes \I_d \bx^{t-1}-\x^{t-1}}^2 \nn, \\
        &\leqtext{\eqref{eq:alg_up_cmpct_eq_2}}  3 \lp \theta^t \rp^2 + 3n \alpha^2 \norm{ \bxh^{t-1} - \bx^{t-1}}^2 + 3 \lp \theta^{t-1} \rp^2 \nn, \\
        &\leqtext{Lemma~\ref{lem:consnss_bound}} 3 \lb (1 + \eta_1) \lamW^2 \lp \theta^{t-1} \rp^2 + \lp 1 +\frac{1}{\eta_1}\rp \alpha^2 \lamW^2 \norm{\delb^{t-1}}^2 \rb \nn \\
        &\qquad + 3 \alpha^2 \norm{ \delb^{t-1}}^2 + 3\lp \theta^{t-1} \rp^2  \nn, \\
        &= 3 \alpha^2 \lp \lp 1 +\frac{1}{\eta_1}\rp \lamW^2 + 1 \rp \norm{ \delb^{t-1}}^2  + 3 ((1 + \eta_1) \lambda^2 +1) \lp \theta^{t-1} \rp^2 \nn. 
    \end{align}
    Finally, we get
    \begin{align}
         \varepsilon^t &\leq 
          \lamW^2 \lp 1 + \gamma_1\rp \varepsilon^{t-1} + \lamW^2 \beta^2 \lp \frac{2}{\gamma_1} + 3 \rp \upsilon^{t-1} + 3 \lamW^2 L^2 ((1 + \eta_1) \lamW^2 +1) \lp \frac{2}{\gamma_1}   + 3 \rp  \EE \lb (\theta^t)^2 \rb    \nn \\
        & + 3 \lamW^2 \beta^2 \bar{\sigma}^2  +3 \alpha^2 \lamW^2 L^2 \lp \frac{2 }{\gamma_1}   + 3\rp \lp \lamW^2 \lp 1 +\frac{1}{\eta_1}\rp + 1\rp \EE \lb \norm{\delb^{t-1}}^2 \rb, \nn 
    \end{align}
    where $\gamma_1, \eta_1 >0$ are Young's parameters. 
    By considering $\gamma_1= \frac{1-\lamW^2}{2 \lamW^2}, \eta_1 = 1$, we get the desired bound.
\end{IEEEproof}

\section{Proof of Lemma~\ref{lem:sum_cnsnss_cmutatv}} \label{pr:lem:sum_cnsnss_cmutatv}
\begin{IEEEproof}
    Substituting $\eta_1 = \frac{1 - \lamW^2}{2 \lamW^2}$ in Lemma~\ref{lem:consnss_bound}, we get
    \begin{align}
       (\theta^t)^2 \leq \lp \frac{1 + \lamW^2 }{2}\rp \lp \theta^{t-1}\rp^2 + \frac{2 \alpha^2 \lamW^2 }{\lp 1  - \lamW^2 \rp} \norm{\delb^{t-1}}^2 \nn.
    \end{align}
   Using Lemma~\ref{lem:seq_bnd}, we have 
    \begin{align}
         \sum_{t=1}^{T}  (\theta^t)^2 &\leq  \frac{2}{1-\lamW^2} (\theta^1)^2 + \frac{4 \alpha^2 \lamW^2 }{\lp 1 -  \lamW^2 \rp^2} \sum_{t=1}^{T-1} \norm{\delb^t}^2 \nn. 
     \end{align}
    From the initialization, we have $\theta^1=0$ and substituting this yields the desired result.
\end{IEEEproof}
    
\section{Proof of Lemma~\ref{lem:HSGD}} \label{pr:lem:HSGD}
\begin{IEEEproof}
     On summing \eqref{eq:lem:3re3} for $1 \leq t \leq T$ and applying \eqref{eq:lem:seq_bnd_frm1} with $\eta_2 = 1$, we obtain,
    \begin{align}
       \sum_{t=1}^T & \phi^t \leq \frac{\phi^{1}}{ 1 - (1-\beta)^2}  +  \frac{6 L^2(1- \beta)^2\alpha^2}{n^2(1 - (1-\beta)^2)} \sum_{t=1}^{T-1}  \EE \lb \norm{\delb^t }^2 \rb  +\frac{12 L^2(1- \beta)^2}{(1 - (1-\beta)^2)n^2} \sum_{t=0}^{T} \EE \lb  \lp \theta^t \rp^2 \rb \nn \\
       &+ \frac{2 \beta^2 \bar{\sigma}^2 T}{n^2(1 - (1-\beta)^2)}. \nn
    \end{align}
     Observing that $\frac{1}{1 -(1-\beta)^2} \leq \frac{1}{\beta}$ for $\beta \in (0,1)$ we have,
     \begin{align}
       \sum_{t=1}^T & \phi^t \leq \frac{\phi^{1}}{ \beta } + \frac{2 \beta \bar{\sigma}^2 T}{n^2} +  \frac{6 L^2 \alpha^2}{n^2 \beta} \sum_{t=1}^{T-1}  \EE \lb \norm{\delb^t }^2 \rb +\frac{12 L^2}{\beta n^2} \sum_{t=0}^{T} \EE \lb  \lp \theta^t \rp^2 \rb. \nn
    \end{align}
    Also, based on the initialization of $\z_i^1$ and Assumption~\eqref{assm:bounded_var}, we have:
    \begin{align}
        \phi^1 &= \EE\lb \norm {\frac{1}{n} \sumin \z_i^1 - \frac{1}{n} \sumin \nabla u_i (\x_i^1)}^2 \rb = \EE\lb \norm {\frac{1}{n} \sumin \frac{1}{b_0} \sum_{r=1}^{b_0} \lp \nfi(\x_i^1, \xi_i^{1,r})  - \nabla u_i (\x_i^1)\rp}^2 \rb,  \nn \\ 
        &\leqtext{(i)}  \frac{\bar{\sigma}^2}{n^2 b_0}. \nn
    \end{align}
    In (i), we have applied Assumption~\ref{assm:bounded_var} and the fact that stochastic local gradient oracles at each node are independent. Substituting $\phi^1$ yields the desired result. The second result can also be obtained by starting from \eqref{eq:lem:3re4} and following similar steps.
\end{IEEEproof}

\section{Proof of Lemma~\ref{lem:sumy}}\label{pr:lem:sumy}
\begin{IEEEproof}
    On summing the bound obtained in  Lemma~\ref{lem:bndnss_y} for $1 \leq t \leq T$ and applying \eqref{eq:lem:seq_bnd_frm1}, we have 
    \begin{align}
        \sum_{t=1}^{T} &\varepsilon^t
        \leq \frac{2}{1- \lamW^2}  \varepsilon^{1} + \frac{8 \beta^2  \lamW^2  }{\lp 1- \lamW^2 \rp^2} \sum_{t=2}^{T} \upsilon^{t-1} + \frac{6 \lamW^2 \beta^2 \bar{\sigma}^2 T}{1- \lamW^2} + \frac{72 \lamW^2 L^2}{\lp 1- \lamW^2 \rp^2}  \sum_{t=2}^{T} \EE \lb (\theta^{t-1})^2 \rb  \nn \\
        &\qquad + \frac{72 \alpha^2 \lamW^2 L^2 }{\lp 1- \lamW^2 \rp^2}  \sum_{t=2}^{T} \EE \lb \norm{ \delb^{t-1}}^2 \rb \nn, \\
        &\leqtext{Lemma~\eqref{lem:HSGD}} \frac{2}{1- \lamW^2}  \varepsilon^{1} + \frac{72 \alpha^2 \lamW^2 L^2 }{\lp 1- \lamW^2 \rp^2}  \sum_{t=2}^{T} \EE \lb \norm{ \delb^{t-1}}^2 \rb  + \frac{72 \lamW^2 L^2}{\lp 1- \lamW^2 \rp^2}  \sum_{t=2}^{T} \EE \lb (\theta^{t-1})^2 \rb  \nn \\
        &\qquad+\frac{8 \beta^2  \lamW^2  }{\lp 1- \lamW^2 \rp^2} \lp \frac{\bar{\sigma}^2}{b_0 \beta}  + 2 \beta \bar{\sigma}^2T    + \frac{6L^2 \alpha^2 }{\beta} \sum_{t=1}^{T-1} \EE \lb \norm{\delb^t}^2  \rb  + \frac{12L^2 }{\beta} \sum_{t=1}^{T} \EE \lb   \lp \theta^t \rp^2 \rb  \rp + \frac{6 \lamW^2 \beta^2 \bar{\sigma}^2 T}{1- \lamW^2}\nn.
    \end{align}
    Further combining the common terms, we get the desired result.
\end{IEEEproof}

\section{Proof of Lemma~\ref{lemU}} \label{pr:lemU}
\begin{IEEEproof}
    Since the surrogate $\tilde{f}$ is a strongly convex function, solving \eqref{eq:alg_up_eq_1} is equivalent to solving a simple convex optimization problem. The optimality condition of convex optimization problem \eqref{eq:alg_up_eq_1} implies:
    \begin{align}
        \ip {\nabla	\tf \lp \xh_i^t, \x_i^t, \xi_i^t \rp + \pi_i^t + \hat{\w}_i^t}{\x_i^t - \xh_i^t} & \geq 0,  \nn
	\end{align}
    where $ \hat{\w}_i^t \in \partial (h + \ind_{\cX}) \mid_{\x = \xh_i^t}$, and $ \ind_{\cX}(\x)$ is an indicator function.
    As per the definition of $\tf \lp \x_i, \x_i^t, \xi_i^t \rp$ given in \eqref{def:tf} and using \eqref{eq:alg_up_eq_3}, we have $\nabla \tf(\x_i^t,\x_i^t,\xib_i^t) = \z_i^{t} $. Furthermore, by adding and subtracting $\nabla \tf(\x_i^t,\x_i^t,\xib_i^t) = \z_i^{t} $ and substituting $\pi_i^t = \y_i^t - \z_i^t $, we get
    \begin{align}
         \ip {\nabla	\tf \lp \xh_i^t, \x_i^t, \xi_i^t \rp + \y_i^t - \z_i^t + \z_i^t -\nabla \tf \lp \x_i^t, \x_i^t, \xi_i^t \rp + \hat{\w}_i^t }{\x_i^t - \xh_i^t} & \geq 0. \nn
    \end{align}
    From the definition of $\tf \lp \x_i, \x_i^t, \xi_i^t \rp$ \eqref{def:tf}, and using \eqref{eq:alg_up_eq_3}, we get $\nabla \tf \lp \hx_i^t, \x_i^t, \xi_i^t \rp - \z_i^t = \mu_i(\xh_i^t - \x_i^t) $, substituting this gives us:
    \begin{align}
         \ip {\mu_i(\xh_i^t - \x_i^t) + \y_i^t+ \hat{\w}_i^t}{ \xh_i^t - \x_i^t} & \leq 0, \nn \\
         \ip { \y_i^t+ \hat{\w}_i^t}{ \xh_i^t - \x_i^t} & \leq - \mu \norm{\xh_i^t - \x_i^t}^2. \nn
    \end{align}
    Further, on summing the above inequality over all $i$ we get
    \begin{align}     
           \frac{1}{n}  \sumin   \ip { \y_i^t}{ \xh_i^t - \x_i^t} +  \frac{1}{n}  \sumin   \ip { \hat{\w}_i^t}{ \xh_i^t - \x_i^t} & \leq  \frac{-1}{n}  \sumin \mu \norm{\xh_i^t - \x_i^t}^2. \label{eq:optimality}
    \end{align}
    From the update equation \eqref{eq:alg_up_eq_2}, convexity of $h +  \ind_{\cX}$ and Assumption~\eqref{assm:condn_on_graph} we obtain
     \begin{align}
       \frac{1}{n} \sumin \lb h(\x_i^{t+1}) + \ind_{\cX}(\x_i^{t+1}) \rb &=  \frac{1}{n} \sumin h \lp \sumjn W_{i,j} \lp \x_j^t + \alpha \lp \xh_j^t - \x_j^t \rp \rp \rp \nn \\
       &\qquad \qquad  +\frac{1}{n} \sumin \ind_{\cX} \lp \sumjn W_{i,j} \lp \x_j^t + \alpha \lp \xh_j^t - \x_j^t \rp \rp \rp, \nn \\ 
    &\overset{(i)}{\leq} \frac{1}{n} \sumin \sumjn W_{i,j}  h \lp  \x_j^t + \alpha \lp \xh_j^t - \x_j^t  \rp \rp +\frac{1}{n} \sumin \sumjn W_{i,j} \ind_{\cX}  \lp \x_j^t + \alpha \lp \xh_j^t - \x_j^t \rp  \rp, \nn \\
    &= \frac{1}{n} \sumjn \lp \sumin W_{i,j} \rp h \lp (1- \alpha) \x_j^t + \alpha \xh_j^t \rp \nn \\
       &\qquad \qquad  +\frac{1}{n} \sumjn\lp \sumin W_{i,j} \rp \ind_{\cX} \lp (1- \alpha) \x_j^t + \alpha \xh_j^t \rp,  \nn \\  
    &\overset{(ii)}{\leq} \frac{1}{n} \sumjn \lp (1- \alpha) h \lp \x_j^t \rp + \alpha h \lp \xh_j^t \rp \rp  +\frac{1}{n} \sumjn \lp (1- \alpha) \ind_{\cX} \lp  \x_j^t \rp + \alpha \ind_{\cX} \lp \xh_j^t \rp \rp. \nn \\
     &= \frac{(1- \alpha) }{n}  \sumjn \lp  h \lp \x_j^t\rp  + \ind_{\cX} \lp  \x_j^t \rp \rp  +\frac{\alpha }{n}  \sumjn \lp h \lp \xh_j^t \rp + \ind_{\cX} \lp \xh_j^t \rp \rp. \nn
    \end{align}
    In $(i)$ and $(ii)$, we have applied the zeroth order convexity condition of $h + \ind_\cX $, property of $\W$ being doubly stochastic, and $W_{i,j}>0$ for all $i,j \in \cV$ \eqref{assm:condn_on_graph}.
    
    From the first order convexity condition of $h(\x) +  \ind_\cX (\x_i^t)$, we have $$h(\xh_i^t) + \ind_\cX (\xh_i^t) \leq h(\x_i^t) +\ind_\cX (\x_i^t) + \ip{ \hat{\w}_i^t }{\xh_i^t - \x_i^t},$$
    Using this we can simplify further as below:
    \begin{align}
       \frac{1}{n} \sumin \lb h(\x_i^{t+1}) + \ind_{\cX}(\x_i^{t+1}) \rb &\leq  \frac{(1- \alpha) }{n}  \sumjn \lp  h \lp \x_j^t\rp + \ind_{\cX} \lp  \x_j^t \rp \rp   +\frac{\alpha }{n}  \sumjn \lp h(\x_j^t) +\ind_\cX (\x_j^t) + \ip{ \hat{\w}_j^t }{\xh_j^t - \x_j^t} \rp, \nn \\
    \frac{1}{n} \sumin \lb h(\x_i^{t+1}) + \ind_{\cX}(\x_i^{t+1}) \rb               &\leq  \frac{1}{n}  \sumjn \lp  h \lp \x_j^t\rp + \ind_{\cX} \lp  \x_j^t \rp \rp  + \frac{\alpha }{n}  \sumjn \ip{ \hat{\w}_j^t }{\xh_j^t - \x_j^t}. \label{eq:convexity_subradient}
    \end{align}
    Dividing \eqref{eq:convexity_subradient} by $\alpha$ and adding in \eqref{eq:optimality} we get
    \begin{align}
     \frac{1}{\alpha n} &\sumin \lb h(\x_i^{t+1}) + \ind_{\cX}(\x_i^{t+1}) \rb +\frac{1}{n}  \sumin   \ip { \y_i^t}{ \xh_i^t - \x_i^t} +  \frac{1}{n}  \sumin   \ip { \hat{\w}_i^t}{ \xh_i^t - \x_i^t} \nn \\
     &\leq  \frac{1}{\alpha n}  \sumin \lp  h \lp \x_i^t\rp + \ind_{\cX} \lp  \x_j^t \rp  \rp  + \frac{1}{n}  \sumjn \ip{ \hat{\w}_j^t }{\xh_i^t - \x_i^t} +\frac{-1}{n}  \sumin \mu \norm{\xh_i^t - \x_i^t}^2. \nn
    \end{align}
    Furthermore, as $\cX$ is a convex set ($g(\x)$ is convex), $\xh_i^t \in \cX$ \eqref{eq:alg_up_eq_1}, Algorithm~\ref{alg:DMSCA} is initialized with a feasible point, and update equation \eqref{eq:alg_up_cmpct_eq_2} is a convex combination of vectors within the set $\cX$. Therefore, we conclude that $\x_i^t \in \cX$ and $\ind_{\cX} \lp \x_i^t \rp = 0$ for all $i \in \cV$ and $t > 0$.  Using this, we obtain:
    \begin{align}
        \frac{1}{n}  \sumin   \ip { \y_i^t}{ \xh_i^t - \x_i^t} 
        \leq  \frac{-1}{n}  \sumin \mu \norm{\xh_i^t - \x_i^t}^2 - \frac{1}{\alpha n} \sumin h(\x_i^{t+1}) + \frac{1}{\alpha n}  \sumin   h \lp \x_i^t\rp.  \label{eq:descent_ineq}
    \end{align}
    From the smoothness of $u(\x)$ (Assumption~\ref{assm:smoothness}), we have
    \begin{align}
        u( \bx^{t+1}) &\leq u(\bx^t) + \ip{ \nabla u( \bx^t)} { \lp \bx^{t+1}-\bx^t \rp} + \frac{L}{2}\norm{ \bx^{t+1}- \bx^t}^2 \nn 
    \end{align}
    Furthermore, pre-multiplying the update equation~\eqref{eq:alg_up_cmpct_eq_2} by $\frac{1}{n}\lp \bo^\mathsf{T} \otimes \I_d\rp$ gives 
    \begin{align}
        \bx^{t+1} =  \bx^t + \alpha \lp \bxh^t - \bx^t \rp \tag{As, $\frac{1}{n}\lp \bo^\mathsf{T} \otimes \I_d\rp \Wu = \frac{1}{n}\lp \bo^\mathsf{T} \otimes \I_d\rp$ } 
    \end{align}
    Using this we can further simplify the quadratic upper bound of $u(\x)$, as below:
    \begin{align}
       u(\bx^{t+1}) &\leq u(\bx^t) + \alpha \Bigg\langle \nabla u( \bx^t) -\y_i^t + \y_i^t, \lp \frac{1}{n} \sumin \xh_i^{t}- \frac{1}{n} \sumin \x_i^t \rp \Bigg\rangle + \frac{\alpha^2 L}{2}\norm{ \bxh^{t}- \bx^t}^2, \nn \\
       &= u(\bx^t) + \frac{\alpha}{n} \sumin \ip{ \nabla u( \bx^t) - \y_i^t} { \lp \xh_i^{t}- \x_i^t \rp} + \frac{\alpha}{n} \sumin \ip{ \y_i^t} { \lp \xh_i^{t}- \x_i^t \rp} +\frac{\alpha^2 L}{2}\norm{ \bxh^{t}- \bx^t}^2, \nn \\
       &\leqtext{\eqref{eq:descent_ineq}} u(\bx^t) + \frac{\alpha}{n} \ip{ (\bo \otimes \I_d) \nabla u( \bx^t) - \y^t} { \xh^{t}- \x^t } +\frac{\alpha^2 L}{2}\norm{ \bxh^{t}- \bx^t}^2\nn \\
       &\qquad + \frac{\alpha}{n} \lp - \sumin \mu \norm{\xh_i^t - \x_i^t}^2 - \frac{1}{\alpha} \sumin h(\x_i^{t+1}) + \frac{1}{\alpha}  \sumin   h \lp \x_i^t\rp \rp, \nn \\
       u(\bx^{t+1}) +\frac{1}{n} \sumin h(\x_i^{t+1}) &- u(\bx^t) - \frac{1}{n}  \sumin   h \lp \x_i^t\rp   \nn \\
       &\leq \frac{\alpha}{n} \ip{ (\bo \otimes \I_d) \nabla u( \bx^t) - \y^t} { \xh^{t}- \x^t } +\frac{\alpha^2 L}{2}\norm{ \bxh^{t}- \bx^t}^2 - \frac{\alpha \mu}{n} \sumin \norm{\xh_i^t - \x_i^t}^2.  \nn 
    \end{align}
    Further applying Cauchy Schwartz inequality and peter-paul's inequality for $\gamma_1>0$ we get,
    \begin{align}
       & u(\bx^{t+1}) +\frac{1}{n} \sumin h(\x_i^{t+1})- u(\bx^t) - \frac{1}{n}  \sumin   h \lp \x_i^t\rp \nn \\
       &\leq \frac{\alpha \gamma_1}{2n} \norm{ (\bo \otimes \I_d) \nabla u( \bx^t) - \y^t}^2 + \frac{\alpha}{2n \gamma_1}  \norm{ \xh^{t}- \x^t }^2 +\frac{\alpha^2 L}{2}\norm{ \bxh^{t}- \bx^t}^2 - \frac{\alpha \mu}{n} \norm{\xh^t - \x^t}^2  \nn \\
       &\overset{\eqref{eq:sub:lem:basic_decentralized_result_4}}{\leq} \frac{\alpha \gamma_1}{2n} \norm{ (\bo \otimes \I_d) \nabla u( \bx^t) -(\bo \otimes \I_d) \bnu(\x^t) + (\bo \otimes \I_d) \bnu (\x^t) - (\bo \otimes \I_d) \by^t + (\bo \otimes \I_d) \by^t - \y^t} \nn \\
       &+ \frac{\alpha}{n} \lp - \mu +\frac{1}{2 \gamma_1 } + \frac{\alpha L}{2} \rp \norm{ \xh^{t}- \x^t } \nn 
    \end{align}
    By further applying the Cauchy-Schwarz inequality and the Peter-Paul inequality for $\gamma_1 > 0$, we obtain:
    \begin{align}
       u( \bx^{t+1}) +\frac{1}{n} \sumin h(\x_i^{t+1}) &- u(\bx^t) - \frac{1}{n}  \sumin   h \lp \x_i^t\rp \nn \\
       &\leq  \frac{\alpha}{2 n}   \gamma_1 \lp 3L^2 \norm{\x^t - \bou \bx^t}^2 +  3 n \norm{ \bnu( \x^t) - \bz^t}^2 +  3\norm{ \y^t - \bou \by^t}^2 \rp \nn \\
        &\quad+ \frac{\alpha}{n} \lp - \mu +\frac{1}{2 \gamma_1 } + \frac{\alpha L}{2} \rp\norm{ \xh^{t}- \x^t }^2 \nn.
    \end{align}
   Taking the Expectation on both sides and summing over $1\leq t \leq T$
    \begin{align}
       &\sumtT u( \bx^{t+1}) +\frac{1}{n} \sumtT \sumin h(\x_i^{t+1}) - \sumtT u(\bx^t) - \frac{1}{n}  \sumtT \sumin   h \lp \x_i^t \rp \nn \\
       &\leq  \frac{\alpha}{2 n}   \gamma_1 \lp 3L^2 \sumtT \EE \lb \norm{\x^t - \bou \bx^t}^2 \rb  +  3 n \sumtT \EE \lb \norm{ \bnu( \x^t) - \bz^t}^2 \rb  +  3 \sumtT \EE \lb \norm{ \y^t - \bou \by^t}^2 \rb \rp \nn \\
        &\quad+ \frac{\alpha}{n} \lp - \mu +\frac{1}{2 \gamma_1 } + \frac{\alpha L}{2} \rp \sumtT \EE \lb \norm{ \xh^{t}- \x^t }^2 \rb \nn.
    \end{align}
    which, on further simplifications, gives 
    \begin{align}
       u&( \bx^{T+1}) +\frac{1}{n} \sumin h(\x_i^{T+1}) -  u(\bx^1) - \frac{1}{n} \sumin   h \lp \x_i^1 \rp \nn \\
       &\leq  \frac{\alpha}{2 n}   \gamma_1 \lp 3L^2 \sumtT \EE \lb \norm{\x^t - \bou \bx^t}^2 \rb  +  3 n \sumtT \EE \lb \norm{ \bnu( \x^t) - \bz^t}^2 \rb  +  3 \sumtT \EE \lb \norm{ \y^t - \bou \by^t}^2 \rb \rp \nn \\
        &\quad+ \frac{\alpha}{n} \lp - \mu +\frac{1}{2 \gamma_1 } + \frac{\alpha L}{2} \rp \sumtT \EE \lb \norm{ \xh^{t}- \x^t }^2 \rb \nn.
    \end{align}
    Further from the zeroth order convexity condition of $h(\x)$ we have $ h( \bx^{T+1}) \leq \frac{1}{n} \sumin h(\x_i^{T+1})$ and from the initialization of D-MSSCA we have $\x_i^1 = \bx^1$ for all $i \in \cV$, using which we get the desired result.

\end{IEEEproof}

\end{document}